\input ./style/arxiv-general.cfg
\documentclass[MSNbibl,number,citesort,seceqn,dvips]{arxbj}
\makeatletter
   \@ifpackageloaded{graphicx}{}{\usepackage{graphicx}}
\makeatother
\usepackage{upgreek}
\usepackage{mathrsfs}

\volume{22}
\issue{3}
\pubyear{2016}
\firstpage{1383}
\lastpage{1430}
\doi{10.3150/15-BEJ696}
\docsubty{FLA}

\makeatletter
\newcommand{\rrvert}{\vert}
\newcommand{\rrVert}{\Vert}
\newcommand{\llvert}{\vert}
\newcommand{\llVert}{\Vert}
\newcommand{\scrB}{\mathscr{B}} 
\newcommand{\scrF}{\mathscr{F}} 
\newcommand{\caB}{\mathcal{B}} 
\newcommand{\caC}{\mathcal{C}} 
\newcommand{\caF}{\mathcal{F}} 
\newcommand{\caG}{\mathcal{G}} 
\newcommand{\caH}{\mathcal{H}} 
\newcommand{\caI}{\mathcal{I}} 
\newcommand{\caK}{\mathcal{K}} 
\newcommand{\caO}{\mathcal{O}} 
\newcommand{\caP}{\mathcal{P}} 
\newcommand{\D}{\mathbb{D}} 
\newcommand{\N}{\mathbb{N}} 
\newcommand{\R}{\mathbb{R}} 
\renewcommand{\L}{\mathbb{L}}
\newcommand{\ii}{\mathrm{i}}
\newcommand{\E}{\mathbb{E}} 
\renewcommand{\P}{\mathbb{P}} 
\newcommand{\intd}{\mathrm{d}} 
\newcommand{\dom}{\operatorname{Dom}}
\newcommand{\tr}{\operatorname{tr}} 
\newcommand{\id}{\mathrm{id}} 
\newcommand{\bigtimes}{\mathop{\,\mbox{\parbox[c][9pt][b]{18pt}{\fontsize{18}{18}\selectfont{$\times$}}}\!\!}}
\newtheorem{lemma}{Lemma}[section]
\newtheorem{theorem}[lemma]{Theorem}
\newtheorem{proposition}[lemma]{Proposition}
\newtheorem{corollary}[lemma]{Corollary}
\newremark{remark}[lemma]{Remark}
\newproclaim{definition}[lemma]{Definition}
\newremark{example}[lemma]{Example}
\newproclaim{assumption}[lemma]{Assumption}
\makeatother

\begin{document}
\begin{frontmatter}

\title{Integration theory for infinite dimensional volatility
modulated Volterra processes}
\runtitle{Integration theory for VMVP}

\begin{aug}
\author[A]{\inits{F.E.}\fnms{Fred Espen}~\snm{Benth}\thanksref{A}\ead[label=e1]{fredb@math.uio.no}}
\and
\author[B]{\inits{A.}\fnms{Andr\'e}~\snm{S\"u\ss{}}\corref{}\thanksref{B}\ead[label=e2]{suess.andre@freenet.de}}
\address[A]{Department of Mathematics, University of Oslo, PO Box 1053
Blindern, N-0316 Oslo, Norway.\\ \printead{e1}}
\address[B]{Facultat de Matem\`{a}tiques, Universitat de Barcelona,
Gran Via de les Corts Catalanes, 585, E-08007 Barcelona, Spain.
\printead{e2}}
\end{aug}

%
\received{\smonth{10} \syear{2013}}
%
\revised{\smonth{7} \syear{2014}}

%
\begin{abstract}
We treat a stochastic integration theory for a class of Hilbert-valued,
volatility-modulated,
conditionally Gaussian Volterra processes. We apply techniques from
Malliavin calculus to define this stochastic integration as a sum of a
Skorohod integral, where the integrand is obtained by applying an
operator to the original integrand, and a correction term involving the
Malliavin derivative of the same altered integrand, integrated against
the Lebesgue measure. The resulting integral satisfies many of the
expected properties of a stochastic integral, including an It\^o
formula. Moreover, we derive an alternative definition using a
random-field approach and relate both concepts. We present examples
related to fundamental solutions to partial differential equations.
\end{abstract}

%
\begin{keyword}
\kwd{Gaussian random fields}
\kwd{Malliavin calculus}
\kwd{stochastic integration}
\kwd{Volterra processes}
\end{keyword}
\end{frontmatter}


\section{Introduction}\label{secintro}
Let throughout this article $0<T<\infty$ be a finite time horizon, fix
$t\in[0,T]$ and let $\caH_1$, $\caH_2$ and $\caH_3$ be three
separable Hilbert spaces. As the main object of investigation in this
article, we introduce the following process:
%
%
\begin{equation}
X(t)= \int_0^t g(t,s)\sigma(s)\delta B(s),
\label{eqdefX}
\end{equation}
where $B$ is a cylindrical Wiener process on $\caH_1$, $\sigma$ is
stochastic process on a time interval $[0,T]$ with values in $L(\caH
_1,\caH_2)$, not necessarily adapted to the Wiener process $B$ and $g$
is a deterministic function depending on two time parameters such that
$g(t,s)\in L(\caH_2,\caH_2)$ for all $0\leq s< t\leq T$. In order for
the stochastic integral to be well-defined, one has to assume that
$g(t,\cdot)\sigma(\cdot)$ is Skorohod integrable on $[0,t]$ so that
$X(t)$ exists as a random element in $\caH_2$. The aim of this article
is to define a stochastic integral with respect to the stochastic
process $X=(X(t))_{t\in[0,T]}$, that is we want to derive an
integration theory for the integral
%
%
\begin{equation}
Z(t) = \int_0^t Y(s) \,\intd X(s),
\label{eqdefZ}
\end{equation}
where we assume that $Y(t)\in L(\caH_2,\caH_3)$ for all $t\in[0,T]$.
With the integration concept we develop, we will see
that naturally $Z(t)\in\caH_3$ for all $t\in(0,T]$.

We want to point out some remarkable facts about the objects we have
just introduced. First, $g(s,s)$ does not have to be defined for any
$s\in[0,T]$, it can be singular on the diagonal. Note, moreover, that
one could put suitable measurability conditions on $\sigma$, such as
predictability, but in general this is not necessary. In such a case,
the integral in (\ref{eqdefX}) would turn out to be an It\^o integral
in a Hilbert space and the condition for the existence of the integral
would be
\[
\E\biggl[\int_0^t \bigl\llVert g(t,s)
\sigma(s)\bigr\rrVert_{L_2(\caH_1,\caH_2)}^2\,\intd s \biggr] < \infty.
\]
At this point, we note that our considerations go beyond the classical
semimartingale case as treated in \cite{protter} for real-valued and
\cite{dapratozabczyk} for Hilbert-valued (semi)martingales. We will
however see in Section~\ref{secSMG} that under some conditions on $g$
(in particular that $g(s,s)$ exists for all $s\in[0,T]$), $X$ will
turn out to be a semimartingale and that the integral in (\ref
{eqdefZ}) and the classical integral with respect to a semimartingale,
denoted by $Y\cdot X$ coincide in some cases. In fact, if $X$ is a
semimartingale, then the difference between these two integrals can be
compared to the situation of the It\^o and Skorohod integral with
respect to the one-dimensional Brownian motion. Next, we list some
examples in order to show the wide range of processes $X$ that can be
used as integrators in (\ref{eqdefZ}).

%
\begin{example}[(Ambit fields)]\label{exambitfields}
The situation that motivates our problem of defining (\ref{eqdefZ})
comes from ambit processes; see \cite{benth}. There we deal with
random fields consisting of a stochastic integral over a random field
defined as follows:
\begin{eqnarray*}
X(t,x) & =& \int_0^t\int_{D(t,x)}
g(t,s,x,y)\sigma(s,y) W(\intd s,\intd y),
\end{eqnarray*}
where $D(t,x)\subseteq{\mathbb{R}^d}$, $W$ is a Gaussian noise white
in time and
(possibly) correlated in space, $g$ is a deterministic function and
$\sigma$ is a random field. These processes are included in our
setting if interpret these equations in the Hilbert space sense of
\cite{dapratozabczyk} where $\caH_2$ is interpreted as $L^2({\mathbb{R}^d})$
and $X(t,\cdot)$ is assumed to be in some $L^2$ with respect to the
spatial parameter. In Section~\ref{secrandomfieldintegral}, we come
back to this example and also derive an integral with respect to a
random field, similar to \cite{walsh}.
\end{example}

%
\begin{example}[(Gaussian processes and VMBV)]
A one-dimensional subclass of $X$ has already been treated in \cite
{alosnualart}, and more generally in \cite{vmlv}. In the former paper,
the authors study integration with respect to $\R$-valued Gaussian
processes where they assumed $\sigma(s)=1$ for all $s\in[0,T]$. In
the latter paper, the authors considered the possibility of a
nontrivial $\sigma$ and they referred to those processes as in (\ref
{eqdefX}) in one dimension as volatility modulated Volterra processes
driven by Gaussian noise (VMBV). In this article, we generalize both
Gaussian and VMBV processes to infinite dimensions. A particular
example is fractional Brownian motion in infinite dimensions. Choose
for all $s\in[0,T]$, $\sigma(s)=Q^{1/2}$ where $Q$ is a nonnegative,
self-adjoint, trace-class operator, let $H\in(0,1)$ and set
\[
g(t,s) = c_H(t-s)^{H-1/2} + c_H \biggl(
\frac{1}{2}-H \biggr)\int_s^t
(u-s)^{H-3/2} \bigl(1-(s/u)^{1/2-H} \bigr)\,\intd u.
\]
Then $X$ is a Hilbert-valued fractional Brownian motion with Hurst
parameter $H$ and covariance operator $Q$.
\end{example}

%
\begin{example}[(Solutions to S(P)DE)]
Another application is stochastic integration with respect to the
solution to a stochastic differential equation in a Hilbert space. This
includes solutions to SPDEs interpreted in the sense of \cite
{dapratozabczyk}. Let, for instance, $X$ be the mild solution to
\[
\intd X(t) = -AX(t) + \sigma\bigl(X(t)\bigr)\,\intd B(t),
\]
with $X_0=0$, where $A$ is an unbounded linear operator, $\sigma$ is a
deterministic function subject to some regularity conditions and $B$ is
again a cylindrical Wiener process on some Hilbert space; see \cite
{dapratozabczyk}, Chapter~6, for a detailed treatment of these
equations. Then if $-A$ generates a strongly continuous semigroup of
linear operators $(g(t))_{t\in[0,T]}$, the mild solution to this
equation is given by the following integral equation:
\[
X(t) = \int_0^t g(t-s)\sigma\bigl(X(s)\bigr)
\,\intd B(s),
\]
which has the form of (\ref{eqdefX}). With the help of the theory we
develop in this article, we are then able to define a stochastic
integral with respect to this solution $X$. A particular example here
are Ornstein--Uhlenbeck processes in infinite dimensions given by the
SDE $\intd X(t) = -AX(t) + F\,\intd B(t)$, where $F$ is a bounded linear
operator. Let the $C_0$-semigroup generated by $-A$ be denoted by
$g(t,s)=\exp(-(t-s)A)$ so that in this case
\[
X(t) = \int_0^t \exp\bigl(-(t-s)A\bigr)F\,\intd
B(s).
\]
\end{example}

The paper is structured in the following way. In Section~\ref
{secpre}, we list some fundamental results and in the subsequent
Section~\ref{secintegral} we give the motivation for the definition
of our integral. Sections~\ref{secrules}~and~\ref{secito}
are dedicated to showing properties of the integral and an It\^o
formula. In the final Sections~\ref{secrandomfieldintegral} and \ref
{secSPDEconncetion}, we provide a random-field integration approach to
the integral and show their equivalence.

Throughout this article, $C$ denotes a positive generic constant, which
may change from line to line without further notice.

\section{Preliminaries}\label{secpre}
\subsection{Vector measures}\label{secvectormeasures}
This subsection deals with the generalization of measures to set
functions taking values in a Banach space, so called \emph{vector
measures}. The case which is most important for us is when the vector
measure is defined on subsets of $\R_+$ taking values in the space of
the bounded linear operators $L(\caH_2,\caH_2)$, where $\caH_2$ is a
separable Hilbert space. After providing the definition and the
important concept of total variation of a vector measure, we list the
most relevant properties here and refer to \cite
{diesteluhl,dunfordschwartz,sergelang} for more details.

%
\begin{definition}
Let $(\caF,\scrF)$ be a measurable space and let $\caB$ be a Banach
space. A set function $\mu\dvtx \caF\rightarrow\caB$ is called a \emph
{finitely additive vector measure}, or in short \emph{vector measure},
if $\mu(F_1\cup F_2)=\mu(F_1)+\mu(F_2)$ for any two disjoint sets
$F_1,F_2\in\scrF$. Moreover if, for any sequence $(F_n)_{n\in\N
}\subseteq\scrF$ of pairwise disjoint subsets of $\caF$, we have
$\mu(\bigcup_{n=1}^\infty F_n) = \sum_{n=1}^\infty\mu(F_n)$, then
$\mu$ is called a \emph{countably additive vector measure}. Note that
the convergence of the sum takes place in the norm topology of $\caB$.

The (\textit{total}) \textit{variation of a vector measure} $\llvert \mu\rrvert
$ is the set
function on $(\caF,\scrF)$ with values in $\R_+\cup\{\infty\}$
defined for all $F\in\scrF$ by
\[
\llvert\mu\rrvert(F):= \sup_{\pi}\sum
_{A\in\pi} \bigl\llVert\mu(A)\bigr\rrVert_\caB,
\]
where the supremum is taken over all partitions $\pi$ of $F$ into a
finite number of pairwise disjoint sets $A\in\caF$.
If $\llvert \mu\rrvert (\caF)<\infty$, then $\mu$ is said to be a
vector measure
of finite variation.
\end{definition}

The total variation of a vector measure is the smallest of all
nonnegative, additive set functions $\lambda$ such that $\llVert \mu
(F)\rrVert
_\caB\leq\lambda(F)$ for all $F\in\scrF$. If a countably additive
vector measure has finite variation, then $\llvert \mu\rrvert $ is
also countably
additive. For vector measures (not necessarily having finite variation)
an integration theory similar to that for $\R_+$-valued measures can
be developed; see \cite{dunfordschwartz}, Section~III.

We briefly list some properties of the Lebesgue--Stieltjes integral
with respect to a Banach-valued function which are important in the
remaining paper.
Let $g$ be a $\caB$-valued function on a finite or infinite interval
of $\R$, which is assumed to have bounded variation. Then we can
define a vector measure $\mu_g$ on all finite subintervals $[a,b]$ by
$\mu_g([a,b]):=g(b)-g(a)$. Then we extend the measure $\mu_g$ onto
the Borel $\sigma$-field of $\R_+$ as its Lebesgue extension. We
denote the extension by $\mu_g$, too. Next, we fix a function $f\dvtx \R
_+\rightarrow L(\caB,\caB_1)$, where $\caB_1$ is another Banach
space. Then, if $f$ is a $\mu_g$-integrable function over an interval
$[a,b]$, we denote the integral with respect to $\mu_g$ throughout
this article by $\int_a^b f(s)g(\intd s)$ instead of $\int_a^b
f(s)\mu_g(\intd s)$. A special choice for $g$ is the identity on $\R
$. In this case, the integral is also known as \emph{Bochner integral}
or \emph{Pettis integral} depending on measurability properties of $f$
and integrability properties of $\llVert f\rrVert $; see \cite
{diesteluhl}, Chapter~II.

Finally, we provide a notion of absolute continuity of the vector
measure $\mu_g$ with respect to the one-dimensional Lebesgue measure
$\lambda$. For this, consider the measure space $([0,T],\scrB
([0,T]),\lambda)$. Let $g\dvtx [0,T]\rightarrow\caB$ be a function of
finite variation and assume that there exists some function $\phi
\dvtx [0,T]\rightarrow\caB$ such that $g(t) = \int_0^t \phi(s)\,\intd s$.
This holds, for instance, when $g$ is Fr\'echet differentiable with
derivative $\phi$. Assume moreover that $f$ is a $\mu_g$-integrable
function taking values in $L(\caB,\caB_1)$, where $\caB_1$ is
another Banach space. Then we have that for every $A\in\scrB([0,T])$
\[
\int_A f(s)g(\intd s) = \int_A
f(s)\mu_g(\intd s) = \int_A f(s)\phi(s)\,\intd
s;
\]
see \cite{dunfordschwartz}, Section~III.11, for more details.
Throughout the rest of this article, we will apply the facts stated in
this subsection to the measure generated by the $L(\caH_2,\caH
_2)$-valued function $g$ used in~(\ref{eqdefX}), always integrating
with respect to the first time argument while leaving the second one fixed.

\subsection{Multidimensional Stieltjes integration}\label{secstieltjes}
In this subsection, we give a quick reminder about an extension of the
integration theory treated in the previous one. We want to define the
Lebesgue--Stieltjes integral with respect to a function that has more
than one argument; see \cite{dunfordschwartz}, Chapter VII. For this,
we need the concept of (locally) bounded variation in the case of
functions with several variables, so-called BV functions. Let for this
$U\subseteq{\mathbb{R}^d}$ be an open subset and let $\caO_c(U)$ be
the set of
all precompact open subsets of $U$. Then a function $g\dvtx {\mathbb
{R}^d}\rightarrow
\R$ is said to be of locally bounded variation in $U$ if $g\in
L_{\mathrm{loc}}^1(U)$ and
\[
V(g,O):= \sup_{\phi}\biggl\{\int_U g(x)
\operatorname{div}\phi(x)\,\intd x; \phi\in\caC^1_c\bigl(O;{
\mathbb{R}^d}\bigr), \llVert\phi\rrVert_{L^\infty(U)}\leq1 \biggr
\}
\]
is finite for all $O\in\caO_c(U)$. Then the set of all functions with
(locally) bounded variation forms a Banach space, which is
nonseparable. For our needs, the most important property of BV
functions is that they are precisely those integrators with respect to
which one can define a Stieltjes integral of all continuous functions.
So one defines an integral like that by starting with the simple
functions $f(x) = 1_{ A }(x) = \prod_{j=1}^d 1_{ [a_i,b_i] }(x_i)$,
that is, indicator functions of sets $A = \bigtimes
_{j=1}^d [a_j,b_j]$, where $a,b\in{\mathbb{R}^d}$ and $a\leq b$
(coordinatewise). For these functions, we define the Stieltjes integral
with respect to $g$ and define the notation $g(A)$ by
%
%
\begin{equation}
g(A):= \int_{\mathbb{R}^d}f(x)g(\intd x) = \sum
_{j=0}^d (-1)^j \mathop{\mathop{\sum
_{x\in{\mathbb{R}^d}}}_{x_i\in\{b_i,a_i\}; i=
1,\ldots,d}}_{\llvert \{
i; x_i=a_i\}\rrvert = j} g(x).
\label{eqfunctionA}
\end{equation}
Note\vspace*{1pt} that the sum inside is a finite sum with at most ${d\choose j}$
summands and it means that we sum over all those $g(x)$ where there are
exactly $j$ arguments that come from the lower point $a$ and the other
ones come from $b$. As an example, one has for $d=1$ the usual result
$g(A)=g(b)-g(a)$, for $d=2$ one has
\[
g(A) = g(b_1,b_2) - g(b_1,a_2) -
g(a_1,b_2) + g(a_1,a_2),
\]
and for $d=3$, (\ref{eqfunctionA}) becomes
\begin{eqnarray*}
g(A) &= & g(b_1,b_2,b_3) -
g(b_1,b_2,a_3) - g(b_1,a_2,b_3)
- g(a_1,b_2,b_3)
\\
&&{} + g(b_1,a_2,a_3) + g(a_1,b_2,a_3)
+ g(a_1,a_2,b_3) - g(a_1,a_2,a_3).
\end{eqnarray*}
These formulas become much simpler if $g$ is the product of $d$
functions with one argument each, that is, $g(x) = \prod_{j=1}^d
g_i(x_i)$. Then (\ref{eqfunctionA}) can be easily seen to reduce to
\[
g(A) = \prod_{j=1}^d \int
_\R1_{ [a_j,b_j] }(x_j)g(\intd x_j)
= \prod_{j=1}^d \bigl(g_j(b_j)-g_j(a_j)
\bigr),
\]
which is what one expects. This yields a measure $\mu_g$ on the Borel
$\sigma$-field on ${\mathbb{R}^d}$ by considering the Lebesgue
extension of $g$
in (\ref{eqfunctionA}). With respect to this measure, one can now
derive an integration theory for real valued functions in $f\in L^p(\mu
_g)$ where $p\in[1,\infty]$, and for such a function we denote the
integral by
\[
\int_O f(x)g(\intd x) = \int_O
f(x)\mu_g(\intd x),
\]
where $O\in\caO_c(U)$. This integration theory will be used in
Section~\ref{secrandomfieldintegral}, where we define a stochastic
integral with respect to a random field, which has $d+1$ variables, $d$
being the spatial dimension and one the temporal dimension.

\subsection{Hilbert-valued Malliavin calculus}\label{secMalliavin}
In this subsection, we provide some ideas and the main results we need
related to Malliavin calculus. It will, however, not be sufficient to
only look at Malliavin calculus for real-valued random variables or
${\mathbb{R}^d}$-valued random vectors as treated extensively in \cite
{nualart}.
Instead we have to deal with random elements taking values in some
separable Hilbert space. Some sources for this are \cite
{carmonatehranchi}, Section~5, or \cite{grorudpardoux} and references
therein. A more general setting is the one treated in \cite
{pronkveraar} where the authors treat Malliavin calculus for random
elements taking values in some UMD Banach space. In this subsection, we
will without further notice identify the dual of a separable Hilbert
space $\caG^*$ with $\caG$.

Let $\caG$ be a separable Hilbert space and let $(W(h),h\in\caG)$ be
an isonormal Gaussian process; see \cite{nualart}, Section~1.1.1, for
some of its properties. Let $(\Omega,\caF,\P)$ be the probability
space induced by the isonormal process. Furthermore, we choose another
separable Hilbert space $\caG_1$, and we consider the class of
Hilbert-valued smooth random elements $F\in L^2(\Omega;\caG_1)$ given
by $F=f(W(h_1),\ldots,W(h_n))$ for $h_1,\ldots,h_n\in\caG$, $n\in
\N$ and $f\dvtx \R^n\rightarrow\caG_1$ which is infinitely Fr\'echet
differentiable with some boundedness condition, polynomially bounded or
bounded. These functions $f$ are dense in $L^2(\Omega;\caG_1)$. For
these random elements $F$, the Malliavin derivative is given by
\[
DF:= \sum_{j=1}^n \frac{\partial f}{\partial x_j}
\bigl(W(h_1),\ldots,W(h_n)\bigr)\otimes h_j.
\]
Consequently, for smooth random elements $F$ we can interpret its
Malliavin derivative as another random element with values in $L_2(\caG
,\caG_1)$, the space of Hilbert--Schmidt operators from $\caG$ to
$\caG_1$, or equivalently in the tensor product $\caG_1\otimes\caG$
since $L_2(\caG,\caG_1)$ is isomorphic to this tensor product.
Throughout this article, we will mainly work with the first approach,
but occasionally use the second one when it is more convenient. We can
also apply projections onto the coordinates of $\caG_1$ leading to
one-dimensional Malliavin calculus. In fact, for some $l\in\caG_1$
\[
D^lF:= \langle DF,l\rangle_{\caG_1} = \sum
_{j=1}^n \biggl\langle\frac{\partial f}{\partial x_j}
\bigl(W(h_1),\ldots,W(h_n)\bigr),l\biggr\rangle
_{\caG_1} h_j.
\]
A special choice for $l$ in the previous equality is an element from a
CONS of $\caG_1$, denoted by $(e_k)_{k\in\N}$ yielding
$D_kF:=D^{e_k}F$. A similar calculation for a CONS in $\caG$ leads to
directional Malliavin derivatives as in the one-dimensional case.

As\vspace*{1pt} in the real-valued Malliavin calculus, the operator $D$ is closable
in $L^2(\Omega;\caG_1)$ and we define the Malliavin derivative of an
element $F\in L^2(\Omega;\caG_1)$ which can be represented as a limit
of a sequence of smooth $\caG_1$-valued random elements $(F_n)_{n\in
\N}$ to be the limit of the Malliavin derivatives of the elements of
the sequence, that is, $DF:= \lim_{n\rightarrow\infty} DF_n$. This
convergence takes place in $L^2(\Omega;L_2(\caG,\caG_1))$. The space
of all such elements will be denoted by $\D^{1,2}(\caG_1)$ and it has
the norm
\[
\llVert F\rrVert^2_{1,2,\caG_1} = \E\bigl[\llVert F\rrVert
_{\caG_1}^2 \bigr] + \E\bigl[\llVert DF\rrVert
_{L_2(\caG,\caG_1)}^2 \bigr],
\]
where the index $\caG_1$ will be dropped if this does not cause any
confusion. One can define, as in the real-valued case, the spaces $\D
^{k,p}(\caG_1)$ for $k\in\N$ and $p\geq1$.

We also need to define a Hilbert-valued equivalent to the divergence
operator $\delta$. This operator $\delta_{\caG_1}\dvtx L^2(\Omega
;L_2(\caG,\caG_1))\rightarrow L^2(\Omega;\caG_1)$ is defined to be
the adjoint of $D$, that is,
\[
\E\bigl[\langle DF,G\rangle_{L_2(\caG,\caG_1)} \bigr] = \E\bigl[\bigl
\langle F,
\delta_{\caG_1}(G)\bigr\rangle_{\caG_1} \bigr],
\]
for $F\in\D^{1,2}(\caG_1)$ and all $G\in L^2(\Omega;L_2(\caG,\caG
_1))$ for which
\[
\bigl\llvert\E\bigl[\langle DF,G\rangle_{L_2(\caG,\caG_1)} \bigr]\bigr
\rrvert
\leq C \bigl(\E\bigl[\llVert F\rrVert^2_{\caG_1} \bigr]
\bigr)^{1/2}.
\]
From now on, we drop the index $\caG_1$ from the divergence operator
if this does not cause any confusion.

Having defined these two operators we will now collect some calculus
rules which we will rely on in the subsequent sections. First, we see
that $D$ and $\delta$ are unbounded linear operators. This implies
that one can pull bounded linear deterministic operators or functionals
in and out of the Malliavin derivative and the divergence operator. As
in the real-valued Malliavin calculus, there is a product and chain
rule for the Malliavin derivative which in the Hilbert-valued case need
some explications. Let $\caG_1$, $\caG_2$ and $\caG_3$ be separable
Hilbert spaces and let $F\in L_2(\caG_1,\caG_2)$ and $G\in L_2(\caG
_2,\caG_3)$ be two random linear operators which are Malliavin
differentiable. Then $GF\in\D^{1,2}(L_2(\caG_1,\caG_3))$ and
%
%
\begin{equation}
D(GF) = (DG)F + GDF, \label{eqproductrule}
\end{equation}
where this equality has to be interpreted as $(D(GF))h = (DG)F(h) +
G(DF)(h)$ for all $h\in\caG_1$. A similar rule applies for
directional Malliavin derivatives. For smooth Hilbert-valued random
elements $F$ and $G$ this is shown in \cite{grorudpardoux}, Lemma~2.1,
and the general case follows by an approximation procedure by
Hilbert-valued smooth random elements. The chain rule in Hilbert-valued
Malliavin calculus is defined for functions $\phi\dvtx \caG_1\rightarrow
\caG_2$ which are either Fr\'echet differentiable or Lipschitz
continuous. Let $F\in\D^{1,2}(\caG_1)$. Then $\phi(F)\in\D
^{1,2}(\caG_2)$ and
%
%
\begin{equation}
D\phi(F) = \phi'(F)DF, \label{eqchainrule}
\end{equation}
where $\phi$ is Fr\'echet differentiable and $\phi'$ denotes the Fr\'
echet derivative of $\phi$. If $\phi$ is only Lipschitz continuous,
then $D\phi(F) = \bar{\phi}DF$ where $\bar{\phi}$ is a random
linear operator from $\caG_1$ to $\caG_2$ whose norm is almost surely
bounded by the smallest Lipschitz constant of $\phi$.

In our setting with a cylindrical Wiener process $B$ on $\caG_0$, we
can make some simplifications. First, we note that in this setting, the
Hilbert space $\caG$ on which the isonormal Gaussian process is
defined is equal to $L^2([0,T];\caG_0)$ and one can reinterpret the
Malliavin derivative $DF$ as some $L_2(\caG_0,\caG_1)$-valued
stochastic process $(D_tF)_{t\in[0,T]}$ on the time interval $[0,T]$
given by
\[
D_tF:= \sum_{j=1}^n
\frac{\partial f}{\partial x_j}\bigl(W(h_1),\ldots,W(h_n)\bigr)\otimes
h_j(t),
\]
where $h_j\in L^2([0,T];\caG_0)$ for all $1\leq j\leq n$. Therefore,
$DF$ actually denotes an equivalence class of functions from $\Omega
\times[0,T]$ with values in $L_2(\caG_0,\caG_1)$, but one can find a
representative such that $D_tF$ is measurable in $\omega$ for all
$t\in[0,T]$ and that $(DF)(\omega)$ is measurable in $t$ for all
$\omega\in\Omega$, which we denote as \emph{the} Malliavin
derivative of $F$. As in the general case, one can define the spaces
$\D^{k,p}(\caG_1)$, but moreover one can also define the spaces $\L
^{k,p}(\caG_1)$ to be $L^p([0,T];\D^{k,p}(\caG_1))$. In the
classical real-valued Malliavin calculus $\L^{1,2}:=\L^{1,2}(\R)$.
If $k=1$, then the norm is given by
%
%
\begin{equation}
\llVert F\rrVert_{\L^{1,p}(\caG_1)}^p = \int_0^T
\E\bigl[\llVert F_t\rrVert_{\caG
_1}^p \bigr]\,\intd t + \int_0^T\int_0^T
\E\bigl[\llVert D_sF_t\rrVert^p_{L_2(\caG
_0,\caG_1)}
\bigr]\,\intd s\,\intd t. \label{eqdefL}
\end{equation}
If $k\geq2$, then iterated Malliavin derivatives and further integrals
are added to this expression.

For the $\caG_1$-valued divergence operator $\delta_{\caG_1}$ this
has the consequence that it reduces to the \mbox{$\caG_1$-}valued Skorohod
integral and for all $G\in L^2([0,T]\times\Omega;L_2(\caG_0,\caG
_1))$ we write $\int_0^T G_s\delta B_s$ instead of $\delta(G)$. If
moreover $G$ is predictable, then this integral turns out to be the
$\caG_1$-valued It\^o integral.

The last issue we focus on in this subsection is the interplay between
the Hilbert-valued Malliavin derivative and Skorohod integral. First,
we have the general commutator relation $D\delta(u) = u + \delta
(Du)$, similar to \cite{nualart}, equation~(1.46). Let
now $u$ be a stochastic
process in $\L^{1,2}(L_2(\caG_0,\caG_1))$ and we assume that for all
$t\in[0,T]$ the process $(D_tu(s))_{s\in[0,T]}$ is Skorohod
integrable and the process $(\int_0^T D_tu(s)\delta B(s))_{t\in
[0,T]}$ has a version which is in $L^2(\Omega\times[0,T];L_2(\caG
_0,\caG_1))$. This condition holds, for instance, if $u$ is twice
Malliavin differentiable. Then $\int_0^T u(s)\delta B(s)\in\D
^{1,2}(\caG_1)$ and for all $t\in[0,T]$
\[
D_t\int_0^T u(s)\delta B(s) =
u(t) + \int_0^T D_tu(s)\delta
B(s).
\]
Finally, we provide a Hilbert-valued integration by parts formula which
is inspired by \cite{carmonatehranchi}, Theorem 5.2. However, we need
it in a slightly more general setting which is why we include a quick
proof here. Before we start, we fix a notation. Let throughout this
article $\tr_{\caG_0}$ denote the trace of a linear operator $A\dvtx \caG
_0\rightarrow L_2(\caG_0,\caG_1)$ taken only over $\caG_0$, that is,
\[
\tr_{\caG_0}(A):= \sum_{k\in\N} \langle
Ae_k,e_k\rangle_{\caG
_0},
\]
where $(e_k)_{k\in\N}$ is a CONS of $\caG_0$. Consequently, the
object $\tr_{\caG_0}(A)$ takes values in $\caG_1$. From the
definition of the Skorohod integral as the adjoint of the Malliavin
derivative one has that for all $u\in L^2(\Omega\times[0,T];L_2(\caG
_0,\caG_1))$ and $A\in\D^{1,2}(L_2(\caG_1,\caG_2))$
%
%
\begin{equation}
\E\biggl[A\int_0^T u(s)\delta B(s) \biggr]
= \E\biggl[\int_0^T \tr_{\caG_1}
\bigl((D_sA)u(s) \bigr)\,\intd s \biggr], \label{eqdualityMDSI}
\end{equation}
where the integrand is $\caG_2$-valued and the integral is understood
as a Bochner integral. Similarly, one could write the trace outside the
integral, which would yield an $L_2(\caG_0,\caG_0\otimes\caG
_2)$-valued integrand and integral. Now we are in the position to
formulate the integration by parts formula.

%
\begin{proposition}\label{propHIbP}
Let $u\in L^2(\Omega\times[0,T];L_2(\caG_0,\caG_1))$ be in the
domain of the Skorohod integral $\delta_{\caG_1}$ and let $A\in\D
^{1,2}(L_2(\caG_1,\caG_2))$. Then $Au\in\dom(\delta_{\caG_2})$ and
%
%
\begin{equation}
\int_0^t Au(s)\delta B(s) = A\int
_0^t u(s)\delta B(s) - \tr_{\caG
_0} \int
_0^t D_s(A)u(s)\,\intd s,
\label{eqHIbP}
\end{equation}
for all $t\in[0,T]$. Note that under the conditions above the
right-hand side of this equality is an element in $L^2(\Omega;\caG_2)$.
\end{proposition}

\begin{pf}
Assume that $\Psi\in\D^{1,2}(L_2(\caG_2,\R))$. We have to show that
%
%
\begin{equation}
\E\biggl[\Psi\int_0^t Au(s)\delta B(s)
\biggr] = \E\biggl[\Psi\biggl(A\int_0^t u(s)
\delta B(s) - \tr_{\caG_0}\int_0^t
D_s(A)u(s)\,\intd s \biggr) \biggr], \label{eqproofIbP}
\end{equation}
for all such $\Psi$. Then a calculation similar to \cite
{carmonatehranchi}, Proposition~5.3, yields
\begin{eqnarray*}
\E\biggl[\Psi\int_0^t Au(s)\delta B(s)
\biggr] & =& \E\biggl[\tr_{\caG_0} \int_0^t
D_s(\Psi)Au(s)\,\intd s \biggr]
\\
& =& \E\biggl[\tr_{\caG_0} \int_0^t
\bigl(D_s(\Psi A) - \Psi D_s A \bigr)u(s)\,\intd s \biggr]
\\
& =& \E\biggl[\Psi\biggl(A\int_0^t u(s)
\delta B(s) - \tr_{\caG
_0}\int_0^t
D_s Au(s)\,\intd s \biggr) \biggr],
\end{eqnarray*}
where we used (\ref{eqdualityMDSI}) in the first equality. This
implies the assertion.
\end{pf}

\section{Stochastic integration}\label{secintegral}
In this section, we provide an exact definition for the stochastic
integral in (\ref{eqdefZ}) with respect to an integrator as in (\ref
{eqdefX}). We are keen on deriving an integration theory that also
covers singular $g$, that is, where $g(t,t)$ is not well-defined. In
order to motivate the definition of the stochastic integral, we provide
a \emph{heuristic} calculation that shows how each term comes into
play. Throughout this and the following sections, we work under the
following assumption which have already been mentioned in Section~\ref
{secintro}.

%
\begin{assumption}\label{assdefintegral}
Fix $T>0$ and let $t\in[0,T]$, and let $B=(B_t)_{t\geq0}$ be a
cylindrical Wiener process on $\caH_1$. Furthermore, let $g(t,s)\in
L(\caH_2,\caH_2)$ (nonrandom) for all $0\leq s < t\leq T$ and let
$(\sigma(t))_{t\geq0}$ be an $L(\caH_1,\caH_2)$-valued stochastic
process such that $g(t,s)\sigma(s)\in L_2(\caH_1,\caH_2)$ for all
$s\in[0,t)$, and $1_{ [0,t] }(\cdot)g(t,\cdot)\sigma(\cdot)\in
\dom(\delta)$ for all $t\in[0,T]$. Assume that for all $s\in[0,t)$
the $L(\caH_2,\caH_2)$-valued vector measure $g(\intd u,s)$ has
bounded variation on $[u,v]$ for all $0\leq s<u<v\leq t$.
\end{assumption}

We note again that unlike in \cite{vmlv} we have not assumed that
$\sigma$ is predictable with respect to~$B$, so that the integral in
(\ref{eqdefX}) is a genuine Skorohod integral. In the following
derivation, we will first assume that $\caH_3=\R$, the general case
will be discussed shortly after. The basic idea for the calculations
that follow is to expand all operators which appear in (\ref{eqdefX})
and (\ref{eqdefZ}) into their coordinates, perform similar
calculations as in \cite{vmlv} and then reassemble the original
operators to get closed-form expressions for Hilbert-valued random
elements. To this end, we fix $(e_k)_{k\in\N}$ to be a CONS of $\caH
_2$ and $(f_l)_{l\in\N}$ to be a CONS of $\caH_1$. Then $(B^l)_{l\in
\N}:= (\langle B,f_l\rangle_{\caH_1})_{l\in\N}$ is a sequence of
independent, one-dimensional Brownian motions. With the help of these
two CONS, we can expand $X^k(t):=\langle X(t),e_k\rangle$ for all
$k\in\N$ in the following way:
%
%
\begin{eqnarray}\label{eqmotivationintegral1}
X^k(t) & =& \biggl\langle\sum_{l\in\N}
\int_0^t g(t,s)\sigma(s) (f_l)
\delta B^l(s),e_k\biggr\rangle_{\caH_2}
\nonumber\\[-8pt]\\[-8pt]\nonumber
& =& \sum_{l\in\N} \int_0^t
\bigl\langle g(t,s),e_k\bigr\rangle_{\caH
_2}\sigma(s)
(f_l)\delta B^l(s),
\end{eqnarray}
where in the last term $\langle g(t,s),e_k\rangle_{\caH_2}$ denotes
the linear functional from $\caH_2$ to $\R$ which is defined by
$\langle g(t,s),e_k\rangle_{\caH_2}(x):= \langle g(t,s)x,e_k\rangle
_{\caH_2}$ for all $x\in\caH_2$. The reason why we introduce this
notation is to perform calculations in a more intuitive manner. In
fact, this is a deterministic linear functional which commutes with the
Skorohod integral and the Malliavin derivative as mentioned in
Section~\ref{secMalliavin}.

Now we are in the position to motivate the definition of the stochastic
integral. In what follows, we first assume that the random integrand
$Y\dvtx [0,T]\rightarrow L(\caH_2,\caH_3)$ is differentiable, but this
assumption will be removed afterward. In the following
$Y(t)(e_k)=Y(t)e_k$ is just applying the linear operator $Y(t)$ to
$e_k\in\caH_2$. In this case, we obtain
%
%
\begin{eqnarray} \label{eqmotivationintegral2}
\hspace*{-12pt}&& \int_0^t Y(s)\,\intd X(s)
\nonumber
\\
&&\quad  = \sum_{k\in\N} \int_0^t
Y(s) (e_k)\,\intd X^k(s)
\nonumber
\\
&&\quad  = \sum_{k\in\N} \biggl(Y(t) (e_k)X^k(t)
- \int_0^t \frac{\partial
Y}{\partial u}(u)
(e_k) X^k(u)\,\intd u \biggr)
\nonumber
\\
&&\quad = \sum_{k\in\N} Y(t) (e_k)X^k(t)
- \sum_{k,l\in\N}\int_0^t
\frac{\partial Y}{\partial u}(u) (e_k) \int_0^u
\bigl\langle g(u,s),e_k\bigr\rangle_{\caH_2} \sigma(s)
(f_l)\,\intd B^l(s)\,\intd u\qquad
\\
&&\quad = \sum_{k\in\N} Y(t) (e_k)X^k(t)
- \sum_{k,l\in\N} \biggl(\int_0^t
\int_0^u \frac{\partial Y}{\partial u}(u)
(e_k)\bigl\langle g(u,s),e_k\bigr\rangle_{\caH_2}
\sigma(s) (f_l)\delta B^l(s)\,\intd u
\nonumber
\\
&&\hspace*{108pt}\qquad{}+ \int_0^t \int
_0^u D_{s,l} \biggl(
\frac{\partial
Y}{\partial u}(u) (e_k) \biggr)\bigl\langle g(u,s),e_k
\bigr\rangle_{\caH_2} \sigma(s) (f_l)\,\intd s\,\intd u \biggr)
\nonumber
\\
&&\quad  = \sum_{k\in\N} Y(t) (e_k)X^k(t)
\nonumber
\\
&&\qquad{} - \sum_{k,l\in\N} \biggl(\int
_0^t \biggl(\int_s^t
\frac{\partial Y}{\partial u}(u) (e_k)\bigl\langle g(u,s),e_k\bigr
\rangle_{\caH
_2} \,\intd u \biggr) \sigma(s) (f_l)\delta
B^l(s)
\nonumber
\\
&&\hspace*{35pt}\qquad{} + \int_0^t
D_{s,l} \biggl(\int_s^t
\frac{\partial
Y}{\partial u}(u) (e_k)\bigl\langle g(u,s),e_k\bigr
\rangle_{\caH_2} \,\intd u \biggr) \sigma(s) (f_l)\,\intd s
\biggr),\nonumber
\end{eqnarray}
where we first expanded the stochastic integral in $\caH_2$ along the
coordinates of integrand and integrator, substituted (\ref
{eqmotivationintegral1}) and did a series expansion in $\caH_1$. Then
we pulled the linear operator $\partial Y(u)/\partial u$ inside the
stochastic integral using Proposition \ref{propHIbP}, used the
stochastic Fubini's theorem; see \cite{nualart}, Exercise 3.2.7,
pulled the deterministic bounded linear operator $\langle
g(t,s),e_k\rangle$ inside the Malliavin derivative and commuted the
Malliavin derivative and the deterministic integral. Using similar
steps and Proposition \ref{propHIbP} again, we calculate the first
term on the right-hand side of the last expression to be equal to
\begin{eqnarray*}
\sum_{k\in\N} Y(t) (e_k)X^k(t)
& =& \sum_{k,l\in\N} Y(t) (e_k)\int
_0^t \bigl\langle g(t,s),e_k\bigr
\rangle_{\caH_2} \sigma(s) (f_l)\delta B^l(s)
\nonumber
\\[-2pt]
& =& \sum_{k,l\in\N} \biggl(\int_0^t
Y(t) (e_k)\bigl\langle g(t,s),e_k\bigr
\rangle_{\caH_2} \sigma(s) (f_l)\delta B^l(s)
\nonumber
\\[-2pt]
&&\hspace*{24pt} {} + \int_0^t
D_{s,l} \bigl(Y(t) (e_k) \bigr)\bigl\langle
g(t,s),e_k\bigr\rangle_{\caH_2} \sigma(s) (f_l)\,\intd s \biggr).
\end{eqnarray*}
Now we substitute this term into (\ref{eqmotivationintegral2}) to obtain
%
%
\begin{eqnarray} \label{eqderivationintegral3}
&& \int_0^t Y(s)\,\intd X(s)
\nonumber
\\[-2pt]
&&\quad = \sum_{k,l\in\N} \biggl(\int_0^t
\biggl(Y(t) (e_k)\bigl\langle g(t,s),e_k\bigr
\rangle_{\caH_2}
\nonumber
\\[-2pt]
&&\hspace*{44pt}\qquad {}- \int_s^t
\frac{\partial Y}{\partial
u}(u) (e_k)\bigl\langle g(u,s),e_k\bigr
\rangle_{\caH_2} \,\intd u \biggr)\sigma(s) (f_l)\delta
B^l(s)
\nonumber
\\[-2pt]
&&\hspace*{25pt}\qquad {}+ \int_0^t
D_{s,l} \biggl(Y(t) (e_k)\bigl\langle g(t,s),e_k
\bigr\rangle_{\caH_2}\nonumber
\\[-2pt]
&&\hspace*{74pt}\qquad {}- \int_s^t
\frac{\partial Y}{\partial
u}(u) (e_k)\bigl\langle g(u,s),e_k\bigr
\rangle_{\caH_2} \,\intd u \biggr)\sigma(s) (f_l)\,\intd s
\biggr)
\\[-2pt]
&&\quad = \sum_{k,l\in\N} \biggl(\int_0^t
\biggl(Y(s) (e_k)\bigl\langle g(s,s),e_k\bigr
\rangle_{\caH_2}
\nonumber
\\[-2pt]
&&\hspace*{45pt}\qquad {}+ \int_s^t Y(u)
(e_k)\biggl\langle\frac{\partial
g}{\partial u}(u,s),e_k\biggr
\rangle_{\caH_2} \,\intd u \biggr)\sigma(s) (f_l)\delta
B^l(s)
\nonumber
\\[-2pt]
&&\hspace*{25pt}\qquad {}+ \int_0^t
D_{s,l} \biggl(Y(s) (e_k)\bigl\langle g(s,s),e_k
\bigr\rangle_{\caH_2}
\nonumber
\\[-2pt]
&&\hspace*{74pt}\qquad {}+ \int_s^t Y(u)
(e_k)\biggl\langle\frac{\partial
g}{\partial u}(u,s),e_k\biggr
\rangle_{\caH_2} \,\intd u \biggr)\sigma(s) (f_l)\,\intd s
\biggr)
\nonumber
\\[-2pt]
&&\quad = \int_0^t \biggl(Y(s)g(s,s) + \int
_s^t Y(u)g(\intd u,s) \biggr)\sigma(s)\delta
B(s)
\nonumber
\\[-2pt]
&&\qquad {}+ \tr_{\caH_1} \int_0^t
D_s \biggl(Y(s)g(s,s) + \int_s^t
Y(u) g(\intd u,s) \biggr)\sigma(s)\,\intd s,\nonumber
\end{eqnarray}
by performing a deterministic integration by parts procedure using the
fact that we can commute Fr\'echet differentiation and the projection
onto the $k$th coordinate since they are bounded linear operators which
commute with the Fr\'echet derivative. Then we summed up over both CONS
$(e_k)_{k\in\N}$ and $(f_l)_{l\in\N}$.

Next, we briefly treat the case of a general separable Hilbert space
$\caH_3$. For this, fix a CONS of $\caH_3$, denoted by $(d_k)_{k\in
\N}$, and use the expansion $Y(t)=\sum_{k\in\N}\langle
Y(t),d_k\rangle d_k$. Then we perform all the calculations above on
each coordinate of $\caH_3$ separately and at the end sum up again to
obtain a closed expression for the integral. Note that for this summing
up to be true, certain summability conditions on the elements in each
coordinate have to be assumed.

These derivations motivate the definition of the following linear
operator for every $h\in L(\caH_2,\caH_3)$:
%
%
\begin{equation}
\caK_g(h) (t,s):= h(s)g(t,s) + \int_s^t
\bigl(h(u)-h(s) \bigr)g(\intd u,s), \label{eqdefKg}
\end{equation}
where the integral is defined as an integral with respect to an $L(\caH
_2,\caH_2)$-valued vector measure whenever it makes sense; note that
$g(\intd u,s)$ has finite variation on all subintervals $[v,t]$ where
$v>s$ by definition. For all such $h$, $\caK_g(h)(t,s)\in L(\caH
_2,\caH_3)$. Sometimes we will call the operator in (\ref{eqdefKg})
the \emph{kernel} associated to $X$. Under the Assumptions \ref
{assdefintegral} and the ones that follow in Definition \ref
{defkernel}, this linear operator is well-defined. We remark that this
operator is the infinite-dimensional analogon to the one which appears
in \cite{vmlv}, which in turn already appeared in~\cite{alosnualart}.
In some special cases, this operator can be written in a different way.
In fact, if $g(s,s)$ is a well-defined linear operator from $\caH_2$
to $\caH_2$, then
%
%
\begin{equation}
\caK_g(h) (t,s) = h(s)g(s,s) + \int_s^t
h(u)g(\intd u,s). \label{eqdefKg2}
\end{equation}
Note that this is the kernel appearing in (\ref
{eqderivationintegral3}), and now it is obvious that (\ref{eqdefKg})
is a generalization of it. If $g(\cdot,s)$ is absolutely continuous
with respect to the one-dimensional Lebesgue measure on $[0,t]$ with
density $\phi(\cdot,s)$, then obviously
\[
\caK_g(h) (t,s) = h(s)g(t,s) + \int_s^t
\bigl(h(u)-h(s) \bigr)\phi(u,s)\,\intd u,
\]
where the integral is understood as a Bochner integral. This situation
applies in particular if $g(\cdot,s)$ is Fr\'echet differentiable. If
$g$ is homogeneous in the its arguments, that is, $g$ depends on $t$
and $s$ only through their difference, the kernel (\ref{eqdefKg}) can
be rewritten as
\[
\caK_g(h) (t,s) = h(s)g(t-s) + \int_0^{t-s}
h(u+s)g(\intd u).
\]
Going back to (\ref{eqderivationintegral3}), we see that we can
define the stochastic integral in (\ref{eqdefZ}) as
%
%
\begin{eqnarray}\label{eqdefintegral}
&&\int_0^t  Y(s)\,\intd X(s)
:= \int_0^t \caK_g(Y) (t,s)
\sigma(s)\delta B(s) + \tr_{\caH
_1}\int_0^t
D_s \bigl(\caK_g(Y) (t,s) \bigr)\sigma(s)\,\intd s.\qquad
\end{eqnarray}

From the definition in the last line, one can see that the integral
$\int_0^tY(s)\,\mathrm{d}X(s)$ does not depend on the particular choices of the
bases $(e_k)_{k\in\N}$ and $(f_l)_{l\in\N}$. In fact, this can be
shown for the stochastic integral as for the usual It\^o integrals in
separable Hilbert spaces, and the trace terms is also independent of
the choice of basis.

Next, we describe the domain of this integral.

%
\begin{definition}\label{defkernel}
Fix $t\geq0$, let $X$ be defined by (\ref{eqdefX}) and assume Assumption \ref{assdefintegral}. We say that a stochastic process $(Y(s))_{s\in
[0,t]}$ belongs to the domain of the stochastic integral with respect
to $X$, if:
\begin{longlist}[(iii)]
\item[(i)] the process $(Y(u)-Y(s))_{u\in(s,t]}$ is
integrable with respect to $g(\intd u,s)$ almost surely,

\item[(ii)] $s\mapsto\caK_g(Y)(t,s)\sigma
(s)1_{ [0,t] }(s)$ is in the domain of the $\caH_3$-valued
divergence operator $\delta B$, and

\item[(iii)] $\caK_g(Y)(t,s)$ is Malliavin
differentiable with respect to $D_s$ for all $s\in[0,t]$ and the \mbox{$\caH
_3$-}valued stochastic process $s \mapsto\tr_{\caH_1} D_s (\caK
_g(Y)(t,s) )\sigma(s)$ is Bochner integrable on $[0,t]$ almost surely.
\end{longlist}
We denote this by $Y\in\caI^X(0,t)$ and the integral $\int_0^t
Y(s)\,\intd X(s)$ is defined by (\ref{eqdefintegral}).
\end{definition}

Now one may ask about the concrete form of the domain of this integral.
So far, we have been unable to derive a characterization of it which is
a similar problem as in anticipating calculus, where the domain of the
divergence operator $\delta$ is not completely characterized. Instead,
one can identify a subset $\L^{1,2}$ in this domain. This is also the
case for the $X$-integral that for some cases, we can identify a subset
of its domain. In fact, if $\sigma$ is assumed to be Malliavin
differentiable then we can define the subset $\caI^X_{1,2}(0,t)$ which
is given by the set of processes for which the seminorm
\[
\llVert Y\rrVert_{\caI^X_{1,2}(0,t)}:= \bigl\llVert\caK_g(Y) (t,
\cdot)\sigma(\cdot)\bigr\rrVert_{\L^{1,2}(L_2(\caH_1,\caH_3))}
\]
is finite. The equivalence classes of this seminorm depend on the exact
shape of $g$ and $\sigma$. Similarly, one can define a semi-inner
product in an obvious way. The set $\caI^X_{1,2}(0,t)$ is then
included in $\caI^X(0,t)$, since for any $Y\in\caI^X_{1,2}(0,t)$ we have
%
%
\begin{eqnarray} \label{eqcontinuityX}
\E\biggl[\biggl\llVert\int_0^t Y(s)\,\intd
X(s)\biggr\rrVert_{\caH_3}^2 \biggr] & \leq& 2\E\biggl[\biggl
\llVert\int_0^t \caK_g(Y) (t,s)
\sigma(s)\delta B(s)\biggr\rrVert_{\caH_3}^2 \biggr]
\nonumber
\\
&& {} + 2T\E\biggl[\int_0^t
\bigl\llVert\tr_{\caH_1} D_s\caK_g(Y) (t,s)
\sigma(s)\bigr\rrVert^2_{\caH_3}\,\intd s \biggr]
\\
&\leq& C_T\bigl\llVert\caK_g(Y) (t,\cdot)\sigma(
\cdot)\bigr\rrVert^2_{\L^{1,2}(L_2(\caH
_1,\caH_3))},\nonumber
\end{eqnarray}
where we have used the continuity of the Skorohod integral on $\L
^{1,2}(\caH_3)$; see \cite{nualart}, equation~(1.47). Note that the second
term in the second line is a part of the Malliavin derivative of the
integrand of the stochastic integral and, therefore, already included
in the norm estimate for the Skorohod integral. This space $\caI
^X_{1,2}(0,t)$ takes the role of $\L^{1,2}$ from classical Malliavin calculus.

\section{Calculus with respect to the integral}
In the first two subsections, we present some general properties of the
stochastic integral defined in the previous section. Afterward, in
Sections~\ref{secOUprocesses} and \ref{secvolterra} we investigate
some particular cases for the integrand.

\subsection{Basic calculus rules}\label{secrules}
At first, we can conclude that the integral defined in the previous
section is linear. This follows immediately from the linearity of the
Malliavin derivative, divergence operator and the Lebesgue(--Stieltjes)
integral. Formally, we have for $Y,Z\in\caI^X(0,t)$ and two constants
$a,b\in\R$ that $aY+bZ\in\caI^X(0,t)$ and
\[
\int_0^t \bigl(aY(s)+bZ(s)\bigr)\,\intd X(s) =
a\int_0^t Y(s)\,\intd X(s) + b\int
_0^t Z(s)\,\intd X(s).
\]
Another immediate property follows from integrating constants. In fact,
if we choose $\caH_2=\caH_3$ and $Y\equiv\id_{\caH_2}$, which is
easily seen to be in $\caI^X(0,t)$ for all $t\geq0$, then $\caK
_g(Y)(t,s) = g(t,s)$ for all $s<t$ and since $g$ is deterministic, we
have $D\caK_g(Y)(t,s)\equiv0$. This implies
%
%
\begin{equation}
\int_0^t\intd X(s) = \int
_0^t \id_{\caH_2}\,\intd X(s) = \int
_0^t g(t,s)\sigma(s)\delta B(s) = X(t).
\label{eqresultdifference15}
\end{equation}
By combining this with the linearity of the integral for the
deterministic integrand $Y(s)=1_{ [u,v] }\id_{\caH_2}$ we obtain
%
%
\begin{equation}
\int_0^t 1_{ [u,v] }(s)\,\intd X(s) =
\int_0^t 1_{ [u,v] }(s)\id
_{\caH_2} \,\intd X(s) = X(v) - X(u), \label{eqresultdifference2}
\end{equation}
which gives us the intuitive property that the integral over an
indicator function is the increment of the integrator. We have
furthermore that if $0<t<T$ and $Y\in\caI^X(0,t)$ then $Y1_{ [0,t]
}\in\caI^X(0,T)$ and
%
%
\begin{equation}
\int_0^T 1_{ [0,t] }(s)Y(s)\,\intd X(s)
= \int_0^t Y(s)\,\intd X(s). \label{eqproperty4}
\end{equation}
This can be seen by splitting the integral $\int_s^T\cdots g(\intd
u,s)$ into two integrals $\int_s^t\cdots g(\intd u,s)$ and $\int
_t^T\cdots g(\intd u,s)$. Note that this equality does not hold if we
only assume that $Y\in\caI^X(0,T)$ because the fact that $X$ is
Skorohod integrable over $[0,T]$ does not in general imply that
$X1_{ [0,t] }$ is Skorohod integrable over $[0,T]$; see \cite{nualart},
Exercise 3.2.1. However, if $\sigma$ is assumed to be
Malliavin differentiable, this does hold. In fact, in this case $\caK
_g(Y)(T,s)\sigma(s)$ is Malliavin differentiable for all $s\in[0,t]$
which implies that $\caK_g(Y)(t,s)\sigma(s)$ is Skorohod integrable
over $[0,t]$ for all $t\in[0,T]$. Combining (\ref{eqproperty4}) with
the linearity of the integral, we have immediately that for $0\leq
u<v\leq t$ and $Y\in\caI^X(0,u)\cap\caI^X(0,v)$, that $Y1_{ [u,v]
}\in\caI^X(0,t)$ and
\[
\int_0^t Y(s)1_{ [u,v] }(s)\,\intd X(s) =
\int_0^v Y(s)\,\intd X(s) - \int
_0^u Y(s)\,\intd X(s).
\]
Using these basic rules, we can derive more interesting properties for
the stochastic integral with respect to $X$.

%
\begin{proposition}\label{proppropdX}
Assume that Assumption \ref{assdefintegral} holds and that $X$ is
defined by (\ref{eqdefX}). Let $t>0$ and assume $Y\in\caI^X(0,t)$.
\begin{longlist}[(iii)]
\item[(i)] Let $Z$ be a random linear operator from $\caH_3$ to another
separable Hilbert space $\caH_4$ which is almost surely bounded. Then
$ZY\in\caI^X(0,t)$ and
\[
\int_0^t ZY(s)\,\intd X(s) = Z\int
_0^t Y(s)\,\intd X(s)\qquad\mbox{almost surely}.
\]
\item[(ii)] The $X$-integral is local, that is, if $Y=0$ on a measurable set
$A\subseteq\Omega$, then
\[
\int_0^t Y(s)\,\intd X(s) = 0 \qquad\mbox{on $A$}.
\]
\item[(iii)] Let $Y$ be a simple process, that is, $Y
= \sum_{j=1}^{n-1} Z_j1_{ (t_j,t_{j+1}] }$ where $Z_j$ is a random
linear operator from $\caH_2$ to $\caH_3$ which is almost surely
bounded for all $j=1,\ldots,n-1$ and $0 \leq t_1 < \cdots< t_n \leq
t$ is a partition of the interval $[0,t]$. Then $Y\in\caI^X(0,t)$ and
\[
\int_0^t Y(s)\,\intd X(s) = \sum
_{j=1}^{n-1} Z_j \bigl(X(t_{j+1})-X(t_j)
\bigr).
\]
\item[(iv)] Let furthermore $\sigma$ be Malliavin differentiable. Then the
$X$-integral is a continuous linear operator from $\caI^X_{1,2}(0,t)$
to $L^2(\Omega;\caH_3)$.
\end{longlist}
\end{proposition}

\begin{pf}
(i)~Note that by definition of $\caK_g$ we have $\caK_g(ZY)(t,s) =
Z\caK_g(Y)(t,s)$ and by the product rule of Malliavin calculus (\ref
{eqproductrule})
\[
D_s \bigl(Z\caK_g(Y) (t,s) \bigr) =
D_s(Z)\caK_g(Y) (t,s) + ZD_s\caK
_g(Y) (t,s).
\]
This and the Hilbert-valued integration by parts formula (\ref
{eqHIbP}) yield the assertion
\begin{eqnarray*}
&& \int_0^t ZY(s)\,\intd X(s)
\\
&&\quad  = \int_0^t Z\caK_g(Y) (t,s)
\sigma(s)\delta B(s)
\\
&&\qquad {} + \tr_{\caH_1} \int_0^t
\bigl(Z D_s\caK_g(Y) (t,s) + D_s(Z)
\caK_g(Y) (t,s) \bigr)\sigma(s)\,\intd s
\\
&&\quad = Z \biggl(\int_0^t \caK_g(Y)
(t,s)\sigma(s)\delta B(s) + \tr_{\caH
_1}\int_0^t
D_s\caK_g(Y) (t,s)\sigma(s)\,\intd s \biggr)
\\
&&\quad = Z\int_0^t Y(s)\,\intd X(s).
\end{eqnarray*}

(ii)~This claim follows immediately from the fact that the Malliavin
derivative, the Skorohod integral and the Lebesgue(--Stieltjes)
integral are local operators.

(iii)~This is a combination of (i) and (\ref{eqresultdifference2}).

(iv)~This is an obvious conclusion from (i) and the continuity in
(\ref{eqcontinuityX}).
\end{pf}

Another property we want to investigate in this subsection are two
projection equalities with respect to the CONS $(e_k)_{k\in\N}$ of
$\caH_2$. The first one is given by the application $\caK_g(Y)(t,s)$
to $e_k$, which yields
%
%
\begin{equation}
\caK_g(Y) (t,s) (e_k) = \caK_{g(e_k)}(Y) (t,s)
\label{eqprojection1}
\end{equation}
which is easy to see from (\ref{eqdefKg}), where $g(e_k)\in\caH_2$
is defined by $g(t,s)e_k$ for all $0\leq s<t\leq T$. The other
projection equality is integration with respect to $X^k=\langle
X,e_k\rangle_{\caH_2}$, which is a real-valued stochastic process.
Applying the definition of the stochastic integral (\ref
{eqdefintegral}) invoking (\ref{eqmotivationintegral1}) with $\caH
_2=\R$, we see
\[
\int_0^t Y(s)\,\intd X^k(s) = \int
_0^t \caK_{\langle g,e_k\rangle
}(Y) (t,s)\delta B(s) +
\tr_{\caH_1} \int_0^t D_s
\caK_{\langle
g,e_k\rangle}(Y) (t,s)\,\intd s,
\]
where $Y$ is a stochastic process taking values in a Hilbert space
$\caH_3$. Note that if $\caH_3=\R$, the function $s\mapsto\caK
_{\langle g,e_k\rangle}(Y)(t,s)$ is real-valued, so we can relate this
to the one-dimensional case treated in \cite{vmlv}. In order to
recover exactly that situation, one has to assume $\caH_1=\R$, too.

\subsection{Semimartingale condition}\label{secSMG}
In this subsection, we focus on stating a condition under which the
stochastic integral process $t\mapsto X(t)$ is a semimartingale. This
condition is a regularity and smoothness assumption on $g$. The
following proposition is inspired by \cite{vmlv}, Proposition 5.

%
\begin{proposition}\label{propsmgcondition}
Let $t>0$ and assume that $g(t,s)$ is well-defined for all $0\leq s\leq
t$. Furthermore, assume that there is a bi-measurable function $\phi
\dvtx [0,T]\rightarrow L(\caH_2,\caH_2)$ such that
%
%
\begin{equation}
g(t,s) = g(s,s) + \int_s^t \phi(v,s)\,\intd v,
\label{eqdifferentiableg}
\end{equation}
for all $0\leq s\leq t$, where this integral is defined in the sense of
Bochner and
\begin{eqnarray*}
\int_0^t \bigl\llVert g(s,s)\bigr\rrVert
_{L(\caH_2,\caH_2)}^2 \,\intd s < \infty\quad\mbox{and}\quad\int
_0^t\int_0^u
\bigl\llVert\phi(u,s)\bigr\rrVert_{L(\caH_2,\caH
_2)}^2 \,\intd s\,\intd u
< \infty.
\end{eqnarray*}
Suppose furthermore that $\sigma$ is adapted to $B$ and pathwise
locally bounded almost surely. Then $X$ defined by (\ref{eqdefX}) is
a semimartingale with decomposition
%
%
\begin{equation}
X(t) = \int_0^t g(s,s)\sigma(s)\,\intd B(s) +
\int_0^t\int_0^s
\phi(s,u)\sigma(u)\,\intd B(u)\,\intd s. \label{eqSMdecomp}
\end{equation}
\end{proposition}

\begin{pf}
Substituting the definition of $g$ and applying a stochastic Fubini
theorem, we obtain
\begin{eqnarray*}
X(t) & =& \int_0^t g(t,s)\sigma(s)\,\intd B(s)
\\
& =& \int_0^t g(s,s)\sigma(s)\,\intd B(s) + \int
_0^t\int_0^s
\phi(s,u)\sigma(s)\,\intd B(u)\,\intd s,
\end{eqnarray*}
which can be seen to be a semimartingale by \cite{metivier}, Remark 26.4.
\end{pf}

From\vspace*{2pt} (\ref{eqdifferentiableg}), we see that $\phi$ can be
interpreted as some sort of derivative of $g$. In fact, if $g(\cdot
,s)$ is Fr\'echet differentiable, then $\phi(t,s) = \frac{\partial
g}{\partial t}(t,s)$. Note that one can also relax the condition on
$\sigma$ to be adapted, leading to so-called Skorohod semimartingales;
see \cite{nualartpardoux}.

In the case when $X$ is a semimartingale, the question arises in which
relationship our integral and the classical integral with respect to
the semimartingale $X $ (see \cite{protter}) stand to each other. The
following proposition provides a partial answer.

%
\begin{proposition}
Assume that Assumption \ref{assdefintegral} holds and that $X$ is
defined by (\ref{eqdefX}). Assume that the semimartingale in
Proposition \ref{propsmgcondition} hold. Suppose that $Y=(Y(s))_{s\in
[0,T]}$ is a predictable stochastic process, which is integrable with
respect to the semimartingale $X$. Assume moreover that either
$D_sY(s)=0$ for almost all $s\in[0,T]$ or $g(s,s)=0$ for almost all
$s\in[0,T]$. Then $Y\in\caI^X(0,t)$ for all $t\in[0,T]$ and
$(Y\cdot X)_t=\int_0^t Y(s)\,\mathrm{d}X(s)$.
\end{proposition}

\begin{pf}
Using the decomposition of $X$ given in (\ref{eqSMdecomp}), the
stochastic integral of $Y$ with respect to $X$, denoted by $Y\cdot X$,
is given by
\begin{eqnarray}
(Y\cdot X)_t & =& \int_0^t
Y(s)g(s,s)\sigma(s)\,\intd B(s) + \int_0^t Y(s)
\int_0^s \phi(s,u)\sigma(u)\,\intd B(u)\,\intd s
\nonumber
\end{eqnarray}
However, by applying integration by parts and the stochastic Fubini
theorem, we get,
%
\begin{eqnarray} \label{eqSMG=dX}
\hspace*{30pt}& =& \int_0^t Y(s)g(s,s)\sigma(s)\,\intd B(s) +
\int_0^t \int_s^t
Y(u)\phi(u,s)\,\intd u\,\sigma(s)\delta B(s)
\nonumber
\\
&& {} + \tr_{\caH_1}\int_0^t
\int_s^t D_s\bigl(Y(u)\bigr)\phi
(u,s)\sigma(s)\,\intd u\,\intd s
\\
& =& \int_0^t Y(s)\,\intd X(s) -
\tr_{\caH_1}\int_0^t D_s
\bigl(Y(s)\bigr)g(s,s)\sigma(s)\,\intd s.\nonumber
\end{eqnarray}
Therefore, these integrals are equal if $D_sY(s)=0$ or $g(s,s)=0$ for
almost all $s\in[0,t]$.
\end{pf}

The additional condition that $D_sY(s)=0$ for almost all $s\in[0,t]$
holds true, for instance, when $Y$ is a simple process as in Proposition
\ref{proppropdX}(iii), where $Z_{t_i}$ is assumed to
be measurable with respect to $(B_t)_{t\in[0,t_i]}$. Similarly, this
also holds when $Y$ is deterministic, or independent of the driving
Wiener process or lags behind the driving Wiener process.

The difference of the integrals comes from the following fact. Starting
from the predictable simple processes for which $D_sY(s)=0$ for almost
all $s\in[0,t]$, we define the semimartingale integral as the closure
of the simple processes under the $L^2([X]^c)$-norm. For the
$X$-integral, this is not what one does. In fact, predictable simple
functions do not even approximate all the processes in the domain of
the $X$-integral with respect to the $L^2([X]^c)$-norm.

%
\begin{example}
For instance, take $\caH_1=\caH_2=\caH_3=\R$, let $Y$ be a Brownian
motion on $[0,T]$ and let $g\equiv1$ and $\sigma\equiv1$. In this
setup, $\caK_g(Y)(t,s)=B(s)$ and $D_s\caK_g(Y)(t,s)=1$. Therefore,
$Y\in\caI^X(0,t)$ for all $t\in[0,T]$ and
%
%
\begin{equation}
\int_0^t Y(s)\,\intd X(s) = \int
_0^t B(s)\,\intd B(s) + \int_0^t
\intd s = \frac{1}{2} \bigl(B(t)^2 + t \bigr). \label{eqexItoXint}
\end{equation}
But the Malliavin derivative of each predictable simple process $Y_n$
in Proposition \ref{proppropdX}(iii) approximating the
Brownian motion is equal to zero and, therefore, does not converge to $1$.
\end{example}

Another set of integrands on which the integrals coincide are
semimartingales of the form
\[
Y(t) = A_t + \int_0^t u(s)\,\intd
B(s),
\]
where $A$ is a process of bounded variation and we assume that
$D_sA(s)=0$ and $D_su(s)=0$ for almost all $s\in[0,T]$.

\subsection{Deterministic integrands and OU processes}\label{secOUprocesses}
In this subsection, we investigate the integral $\int_0^t h(t,s)\,\intd
X(s)$ if $h$ is chosen to be deterministic. To this end, fix $t>0$ and
let $s\mapsto h(t,s)$ be a measurable function from $[0,t]$ into
$L(\caH_2,\caH_3)$ such that $h(t,u)-h(t,s)$ is integrable with
respect to $g(\intd u,s)$ over $[s,t]$ and Definition \ref
{defkernel}(ii) holds. In this case, the kernel
$\caK_g(h)$ is deterministic and has therefore a Malliavin derivative
equal to zero, which immediately implies condition (iii) in Definition \ref{defkernel}. Then, by the
definition of the integral we have
\begin{eqnarray*}
\int_0^t h(t,s)\,\intd X(s)
&=& \int_0^t \caK_g\bigl(h(t,
\cdot)\bigr) (t,s)\sigma(s)\delta B(s)
\\
& =& \int_0^t \biggl(h(t,s)g(t,s) + \int
_s^t \bigl(h(t,u)-h(t,s) \bigr)g(\intd u,s)
\biggr)\sigma(s)\delta B(s).
\end{eqnarray*}
Note that this is a new process of the same type as $X$ in (\ref
{eqdefX}) where now $\caK_g(h)$ plays the role of~$g$.

Now we apply this to a specific deterministic integrand which is
stationary in time in order to show that Ornstein--Uhlenbeck (OU)
processes driven by $X$ are representable as a stochastic integral of a
deterministic integrand with respect to $X$. By an OU process, we mean
the solution $(Y(t))_{t\geq0}$ to the infinite-dimensional SDE $\mathrm{d}Y(t)
= -AY(t)\,\intd t + F\,\intd X(t)$, or in integral terms
%
%
\begin{equation}
Y(t) = -\int_0^t AY(s)\,\intd s + FX(t),
\label{eqOU}
\end{equation}
where $(X(t))_{t\geq0}$ is the process defined in (\ref{eqdefX}) and
$-A$ is an unbounded linear operator from $\caH_3$ to $\caH_3$ whose
domain is dense in $\caH_3$ and that generates a $C_0$-semigroup on
$\caH_3$ denoted by $S(t)=\mathrm{e}^{-tA}$ for all $t\geq0$.
Moreover, $F$ is a bounded linear operator from $\caH_2$ to $\caH_3$.
Note that this is a first step into the investigation of SDEs driven by
a Volterra-type process. As suggested from the classical situation
where the driving noise in (\ref{eqOU}) is a Wiener process, we
obtain that the solution to (\ref{eqOU}) is given by a stochastic convolution.

%
\begin{proposition}\label{propOU}
Suppose that Assumption \ref{assdefintegral} holds. Assume
furthermore that $t\mapsto\sigma(t)$ is pathwise locally bounded
almost surely. Assume for all $s\in[0,T]$ that the map $(s,T)\ni
u\mapsto(\mathrm{e}^{-(u-s)A}-\mathrm{e}^{-sA})$ is integrable with
respect to $g(\intd u,s)$, that
%
%
\begin{equation}
\E\biggl[\int_0^t\biggl\llVert\int
_0^s F\caK_g(S) (s,u)\sigma(u)
\delta B(u)\biggr\rrVert_{\caH_3}^2\,\intd s \biggr] < \infty,
\label{eqOUcond2}
\end{equation}
and
%
%
\begin{equation}
\lim_{t\downarrow s} \bigl\llVert\bigl(\id_{\caH_2}-\mathrm
{e}^{-(t-s)A} \bigr)g(t,s)\bigr\rrVert_{L_2(\caH_2,\caH_3)} = 0. \label
{eqOUcond}
\end{equation}
Then $S\in\caI^X(0,t)$ and
%
%
\begin{equation}
Y(t) = \int_0^t S(t-s)F\,\intd X(s) = \int
_0^t \mathrm{e}^{-(t-s)A}F\,\intd X(s)
\label{eqsolOU}
\end{equation}
solves the SDE (\ref{eqOU}).
\end{proposition}

First, let us explain the two conditions (\ref{eqOUcond2}) and (\ref
{eqOUcond}). The former one is the Hilbert valued equivalent of the one
in \cite{nualart}, Exercise 3.2.7, which is sufficient for an
application of the stochastic Fubini theorem with the divergence
operator that we will use in the following proof. If $\sigma$ is
assumed to be predictable, the Skorohod integral becomes an It\^o
integral since $\caK_g(h)$ is deterministic. In this case, one can
derive similarly as in \cite{vmlv}, Lemma~2, that
\[
\int_0^T\int_s^T
\bigl\llVert g(u,s)\bigr\rrVert_{L(\caH_2,\caH_2)}^2\,\intd u\,\intd s <
\infty
\]
is a sufficient condition for the stochastic Fubini theorem to hold.
The latter condition (\ref{eqOUcond}) can be interpreted as a kind of
weighted $C_0$-semigroup condition. In fact, we demand that the
semigroup $S(t)=\mathrm{e}^{-tA}$ goes to the identity as $t\downarrow
0$ in a family of norms with weights given by the functions $g(\cdot,s)$.

\begin{pf*}{Proof of Proposition \ref{propOU}}
Let throughout the proof $\zeta\in\dom(A^*)$. By substituting (\ref
{eqsolOU}) and the definition of the stochastic integral into (\ref
{eqOU}) and applying the stochastic Fubini theorem, one obtains
\begin{eqnarray*}
\biggl\langle\zeta,\int_0^tAY(u)\,\intd u\biggr
\rangle_{\caH_3} & =& \int_0^t \biggl
\langle A^*\zeta,\int_0^u
\caK_g(SF) (u,s)\sigma(s)\delta B(s)\biggr\rangle_{\caH_2}\,
\intd u
\\
& =& \int_0^t \biggl\langle\biggl(\int
_s^t A\caK_g(SF) (u,s)\,\intd u
\biggr)^*\zeta,\sigma(s)\delta B(s)\biggr\rangle_{\caH_2},
\end{eqnarray*}
and similarly for the first term of (\ref{eqOU})
\[
\bigl\langle\zeta,Y(t)\bigr\rangle_{\caH_3} = \int_0^t
\bigl\langle\caK_g(SF) (t,s)^*\zeta,\sigma(s)\delta B(s)\bigr
\rangle_{\caH_2}.
\]
These two equalities yield that
\begin{eqnarray}
\label{eqproofOU1} && \biggl\langle\zeta,Y(t) + \int_0^t
AY(u)\,\intd u\biggr\rangle_{\caH_3}
\nonumber\\[-8pt]\\[-8pt]\nonumber
&&\quad = \int_0^t \biggl\langle\biggl(
\caK_g(SF) (t,s)+\int_s^t A\caK
_g(SF) (u,s)\,\intd u \biggr)^*\zeta,\sigma(s)\delta B(s)\biggr\rangle
_{\caH_2}.
\nonumber
\end{eqnarray}
Now we investigate the term in the brackets in the previous equality.
Written out, this term is equal to
%
%
\begin{eqnarray}\label{eqproofOU2}
&& \mathrm{e}^{-(t-s)A}Fg(t,s) + \int_s^t
\bigl(\mathrm{e}^{-(t-u)A}-\mathrm{e}^{-(t-s)A} \bigr)Fg(\intd u,s)
\nonumber\\[-8pt]\\[-8pt]\nonumber
&&\quad{} + \int_s^t \biggl(A\mathrm{e}^{-(u-s)A}Fg(u,s)
+ \int_s^u \bigl(A\mathrm{e}^{-(u-v)A}-A
\mathrm{e}^{-(u-s)A} \bigr)Fg(\intd v,s) \biggr)\,\intd u.
\end{eqnarray}
Using Fubini's theorem and \cite{engelnagel}, Lemma~II.1.3(ii), we see that
%
%
\begin{eqnarray} \label{eqproofOU3}
&& \int_s^t\int_s^u
\bigl(A\mathrm{e}^{-(u-v)A}-A\mathrm{e}^{-(u-s)A} \bigr)Fg(\intd v,s)\,\intd u
\nonumber
\\
&&\quad = \int_s^u\int_v^t
\biggl(-\frac{\partial}{\partial u}\mathrm{e}^{-(u-v)A} + \frac
{\partial}{\partial u}
\mathrm{e}^{-(u-s)A} \biggr)\,\intd uFg(\intd v,s)
\\
&&\quad = -\int_s^t \bigl(\mathrm{e}^{-(t-v)A}-
\mathrm{e}^{-(t-s)A} \bigr)Fg(\intd v,s)+\int_s^t
\bigl(1-\mathrm{e}^{-(v-s)A} \bigr)Fg(\intd v,s).\nonumber
\end{eqnarray}
Note that by an integration by parts procedure and using \cite
{engelnagel}, Lemma~II.1.3(ii), again, the second term on the
right-hand side
of the last equality may be calculated to be
%
%
\begin{eqnarray}\label{eqproofOU4}
&&\int_s^t  \bigl(1-\mathrm{e}^{-(v-s)A}
\bigr)Fg(\intd v,s)
\nonumber\\[-8pt]\\[-8pt]\nonumber
&&\quad  = \bigl(1-\mathrm{e}^{-(t-s)A} \bigr)Fg(t,s) - \int_s^t
A\mathrm{e}^{-(v-s)A}Fg(v,s)\,\intd s,
\end{eqnarray}
where we used the condition (\ref{eqOUcond}) so that the boundary
term at $s$ does not appear. Assembling (\ref{eqproofOU2})--(\ref
{eqproofOU4}), we obtain
\[
\caK_g(SF) (t,s) + \int_u^tA
\caK_g(SF) (s,u)\,\intd s = Fg(t,s),
\]
which in turn, by substituting this into (\ref{eqproofOU1}) and using
the fact that $F$ is bounded, implies the assertion
\begin{eqnarray*}
\biggl\langle\zeta,Y(t) + \int_0^t AY(u)\,\intd
u\biggr\rangle_{\caH_3} & =& \biggl\langle\zeta,\int_0^t
Fg(t,s)\sigma(s)\delta B(s) \biggr\rangle_{\caH_3}
 = \bigl\langle\zeta,FX(t)\bigr\rangle_{\caH_3}.
\end{eqnarray*}

The proof ends by the fact that $A$ and hence $A^*$ are densely defined
in $\caH_3$ which implies~(\ref{eqsolOU})
\end{pf*}

Note that (\ref{eqsolOU}) has the exact same form as in the classical
case when the driving process is a Wiener process. We however have to
remark that this does not show uniqueness of the solution to the SDE
(\ref{eqOU}).

\subsection{Volterra processes as integrands}\label{secvolterra}
In this subsection, we want to calculate the stochastic integral with
respect to $X$ when the integrand is a Volterra process itself. This
will also lead us to an integral of type $\int_0^t X(s)\,\intd X(s)$.
Throughout this subsection, we assume that $\sigma\equiv\id_{\caH
_2}$ and assume that $B$ is a (cylindrical) Wiener process on $\caH
_2$. However, one can reintroduce $\sigma$ into the calculations by
formally considering $\sigma\,\intd B$ as the noise. Moreover, we assume
in the following without loss of generality that $\caH_3=\R$. The
case of general separable Hilbert spaces then follows by the expansion
$Y(s)=\sum_k \langle Y(s),f_k\rangle f_k$, where $(f_k)_{k\in\N}$ is
a CONS of $\caH_3$, and performing the subsequent calculations coordinatewise.

Let in the following $h$ be an $L(\caH_2,L(\caH_2,\R))$-valued
deterministic function in $\caI^X(0,t)$ and by the results of the
previous subsection
%
%
\begin{equation}
Y(s) = \int_0^s h(u)\,\intd X(u) = \int
_0^s \caK_g(h) (s,v)\,\intd B(v).
\label{eqintegrandVolterra}
\end{equation}
Here, the integral is in fact a classical It\^o integral. Furthermore,
let $(e_k)_{k\in\N}$ be a CONS of $\caH_2$.

%
\begin{proposition}\label{prop441}
Suppose that Assumption \ref{assdefintegral} holds and assume that
$\sigma,B,h,Y$ and $\caH_3$ are as stated above. Assume furthermore
that for almost all $s\in[0,T]$, $g(s,s)$ is well-defined as a linear
operator on $\caH_2$. Then $Y\in\caI^X(0,t)$ for all $t\in[0,T]$
and $\int_0^t Y(s)\,\mathrm{d}X(s)$ is a sum of an element in the second and
zeroth Wiener chaos.
\end{proposition}

\begin{pf}
Using the definition of the stochastic integral with respect to $X$ and
the notion of the projection in (\ref{eqprojection1}), we have
%
%
\begin{eqnarray}\label{eqproofvolterra4}
&& \int_0^t Y(s)\,\intd X(s)
\nonumber\\[-8pt]\\[-8pt]\nonumber
&&\quad = \sum_{m\in\N} \biggl(\int_0^t
\caK_g(Y) (t,s) (e_m)\delta B^m(s) + \int
_0^t D_{s,m}\caK_g(Y)
(t,s) (e_m)\,\intd s \biggr).
\end{eqnarray}
Now, taking advantage of the specific form of $Y$, we calculate the
integrand in the stochastic integral in (\ref{eqproofvolterra4}) to be
%
%
\begin{eqnarray}\label{eqproofvolterra2}
&& \caK_{g(e_m)}(Y) (t,s)
\nonumber
\\
&&\quad = Y(s)g(t,s) (e_m) + \int_s^t
\bigl(Y(u)-Y(s) \bigr)g(\intd u,s) (e_m)
\nonumber
\\
&&\quad = \int_0^s \caK_g(h) (t,s)\,\intd B(v)g(t,s) (e_m)
\nonumber
\\
&&\qquad{} + \int_s^t\int
_0^s \bigl(\caK_g(h) (u,v) - \caK
_g(h) (s,v) \bigr)\,\intd B(v)g(\intd u,s) (e_m)
\nonumber
\\
&&\qquad {} + \int_s^t\int
_s^u \caK_g(h) (u,v)\,\intd B(v)g(
\intd u,s) (e_m)\nonumber
\\
&&\quad = \sum_{l\in\N} \biggl(\int_0^s
\caK_{g(e_l)}(h) (t,s)g(t,s) (e_m)\,\intd B^l(v)
\nonumber
\\
&&\hspace*{20pt}\qquad{} + \int_s^t\int
_v^t \caK_{g(e_l)}(h) (u,v)g(\intd u,s)
(e_m)\,\intd B^l(v)
\nonumber
\\
&&\hspace*{20pt}\qquad {} + \int_0^s\int
_s^t \bigl(\caK_{g(e_l)}(h) (u,v) -
\caK_{g(e_l)}(h) (s,v) \bigr)g(\intd u,s) (e_m)\,\intd
B^l(v) \biggr)
\\
&&\quad = \sum_{l\in\N} \biggl(\int_0^s
\caK_{g(e_m)}\bigl(\caK_{g(e_l)}(h) (\cdot,v)\bigr) (t,s)\,\intd
B^l(v)
\nonumber
\\
&&\hspace*{19pt}\qquad {} + \int_s^t\int
_v^t \caK_{g(e_l)}(h) (u,v)g(\intd u,s)
(e_m)\,\intd B^l(v) \biggr)
\nonumber
\\
&&\quad = \sum_{l\in\N} \int_0^t
\tilde{\caK}^{l,m}_{g,g}(h) (s,v,t)\,\intd B^l(v),\nonumber
\end{eqnarray}
where we have set
%
%
\begin{eqnarray}\label{eqproofvolterra5}
\tilde{\caK}^{l,m}_{g_1,g_2}(h) (s,v,t) &= & 1_{ \{v\leq s\} }\caK
_{g_2(e_m)}\bigl(\caK_{g_1(e_l)}(h) (\cdot,v)\bigr) (t,s)
\nonumber\\[-8pt]\\[-8pt]\nonumber
&&{} + 1_{ \{v>s\} }\int_v^t
\caK_{g_1(e_l)}(h) (u,v)g_2(\intd u,s) (e_m).
\end{eqnarray}
Applying the Malliavin derivative operator in (\ref
{eqproofvolterra2}) yields
%
%
\begin{eqnarray}\label{eqproofvolterra3}
D_{s,m}\caK_{g(e_m)}(Y) (t,s) & =& D_{s,m} \biggl(\sum
_{l\in\N}\int_0^t
\tilde{\caK}^{l,m}_{g,g}(h) (s,v,t)\,\intd B^l(v)
\biggr)
\nonumber\\[-8pt]\\[-8pt]\nonumber
& =& \tilde{\caK}^{m,m}_{g,g}(h) (s,s,t).
\end{eqnarray}
Note that on the right-hand side the argument $(s,s,t)$ appears which
implies the explicit appearance of $g(s,s)$ in the formulas; see below
for a comment. Assembling all the terms (\ref
{eqproofvolterra4})--(\ref{eqproofvolterra3}), we obtain
\begin{eqnarray*}
\int_0^t Y(s)\,\intd X(s)
&=& \sum_{m\in\N} \biggl(\int
_0^t \sum_{l\in\N}
\int_0^t \tilde{\caK}^{l,m}_{g,g}(h)
(s,v,t)\,\intd B^l(v)\delta B^m(s) + \int
_0^t \tilde{\caK}^{l,m}_{g,g}(h)
(s,s,t)\,\intd s \biggr)
\\
&= &\sum_{l,m\in\N}\int_0^t
\int_0^t \tilde{\caK}^{l,m}_{g,g}(h)
(s,v,t)\,\intd B^l(v)\delta B^m(s) + \sum
_{m\in\N}\int_0^t \tilde{
\caK}^{m,m}_{g,g}(h) (s,s,t)\,\intd s
\\
&= &\sum_{l,m\in\N}\int_0^t
\int_0^s \bigl(\tilde{\caK
}^{l,m}_{g,g}(h) (s,v,t) + \tilde{\caK}^{m,l}_{g,g}(h)
(v,s,t) \bigr)\,\intd B^l(v)\,\intd B^m(s)
\nonumber
\\
&&{} + \sum_{m\in\N}\int
_0^t \tilde{\caK}^{m,m}_{g,g}(h)
(s,s,t)\,\intd s.
\end{eqnarray*}
If we now define the linear operator $\tilde{\caK}_{g,g}(h)$ by
$\tilde{\caK}_{g,g}(h)(s,v,t)(e_l\otimes e_m):= \tilde{\caK
}^{l,m}_{g,g}(h)(s,v,t)$, we see in total that
%
%
\begin{eqnarray}\label{eqcalcVolterra}
\int_0^t Y(s)\,\intd X(s) &= & \int
_0^t \int_0^s
\bigl(\tilde{\caK}_{g,g}(h) (s,v,t)+\tilde{\caK}_{g,g}(h)
(v,s,t) \bigr)\,\intd B(v)\otimes\intd B(s)
\nonumber\\[-8pt]\\[-8pt]\nonumber
&&{} + \tr_{\caH_2}\int_0^t \tilde{
\caK}_{g,g}(h) (s,s,t)\,\intd s.
\end{eqnarray}
From this identity, we can see that every integrator in this formula is
deterministic and, therefore, we can conclude that the integral $\int
_0^t Y(s)\,\intd X(s)$ belongs to the second Wiener chaos plus a term
from the zeroth Wiener chaos.
\end{pf}

Note that the term $\tilde{\caK}_{g,g}(h)(s,s,t)$ appears in the
correction term in the formula (\ref{eqcalcVolterra}). This implies
that the integral term in (\ref{eqproofvolterra5}) does not appear,
but in the first summand in (\ref{eqproofvolterra5}) the terms
$g(s,s)$ will appear explicitly. This is why we have to assume
sufficient regularity for the kernel $g$, in particular that $g(t,t)$
exists as a linear operator for all $t\in[0,T]$. The reason for $g$ to
be regular is that it appears in the integrand (\ref
{eqintegrandVolterra}) as well as in the integrator. Otherwise, both
might have a singularity at the same point which would cause problems
when integrating.

%
\begin{remark}
Similarly to these calculations above, one can easily extend the
results to iterated $X$ integrals of higher order. As long as $\sigma$
is deterministic, the $X$-integral of an integrand $Y$, which is itself
an iterated $X$-integral, consists of a Skorohod integral over $\caK
_g(Y)$ with respect to a cylindrical Wiener process, and a pathwise
integral of the Malliavin derivative of $\caK_g(Y)$. This increases
the number of the Wiener chaos of the integrand by one (Wiener
integral) and decreases the number of the Wiener chaos of the integrand
by one (pathwise integral). Therefore, for $k$ odd (even), the $k$th
iterated $X$-integral of a deterministic function is an element of all
odd (even and zeroth) Wiener chaoses up to order $k$, respectively.
\end{remark}

Now we apply this to the problem of calculating $\int_0^t X(s)\,\intd
X(s)$. Since both elements are $\caH_2$-valued, we cannot apply (\ref
{eqcalcVolterra}) verbatim. However, by the Riesz representation
theorem, we interpret the integrand as a linear functional on $\caH_2$
and, therefore,
\begin{eqnarray*}
\int_0^t X(s)\,\intd X(s)&:=& \int
_0^t X(s)^*\,\intd X(s)  = \int
_0^t \bigl\langle X(s),\intd X(s)\bigr\rangle
= \sum_{k\in\N} \int_0^t
X^k(s)\,\intd X^k(s).
\end{eqnarray*}

%
\begin{proposition}
Suppose that Assumption \ref{assdefintegral} holds and assume that
$\sigma,B$ and $\caH_3$ are as stated at the beginning of this
section. Assume furthermore that for almost all $s\in[0,T]$, $g(s,s)$
is well-defined as a linear operator on $\caH_2$. Then $X\in\caI
^X(0,t)$ for all $t\in[0,T]$. Moreover, this integral can be written
as the sum of $\frac{1}{2}\langle X(t),X(t)\rangle_{\caH_2}$ and
some correction term.
\end{proposition}

\begin{pf}
In this case, we choose $h(s)=\id_{\caH_2}$ which implies $\caK
_g(h)(t,s) = g(t,s)$. Since $X^k(s) = \int_0^s \langle
g(s,u),e_k\rangle_{\caH_2}\,\intd B(s)$, we obtain that $D_{s,l}X^k(u)
= \langle g(u,s)e_l,e_k\rangle_{\caH_2}$. Then
%
%
\begin{eqnarray} \label{eqproofXX}
&& \sum_{k\in\N}\int_0^t
X^k(s)\,\intd X^k(s)
\nonumber
\\
&&\quad = \sum_{k\in\N}\int_0^t
\biggl(X^k(s)\bigl\langle g(s,s),e_k\bigr\rangle
_{\caH_2} + \int_s^t
X^k(r)\bigl\langle g(\intd r,s),e_k\bigr
\rangle_{\caH
_2} \biggr)\,\intd B(s)
\nonumber
\\
&&\qquad {} + \sum_{k,l\in\N} \int
_0^tD_{s,l} \biggl(X^k(s)
\bigl\langle g(s,s),e_k\bigr\rangle_{\caH_2} + \int
_s^t X^k(r)\bigl\langle g(\intd
r,s),e_k\bigr\rangle_{\caH_2} \biggr) (e_l)\,\intd s\qquad
\\
&&\quad = \int_0^t \int_0^s
\bigl\langle g(s,u)\,\intd B(u),g(s,s)\,\intd B(s)\bigr\rangle_{\caH_2}
\nonumber
\\
&&\qquad {}+ \int_0^t\int
_s^t\biggl\langle\int_0^r
g(r,u)\,\intd B(u),g(\intd r,s)\,\intd B(s)\biggr\rangle_{\caH_2}
\nonumber
\\
&&\qquad {} + \tr_{\caH_2}\int_0^t
\biggl(\bigl\langle g(s,s),g(s,s)\bigr\rangle_{\caH_2} + \int
_s^t \bigl\langle g(r,s),g(\intd r,s)\bigr
\rangle_{\caH_2} \biggr)\,\intd s.\nonumber
\end{eqnarray}
We want to find a formula that links this term with $\langle
X(t),X(t)\rangle_{\caH_2}$, which can be calculated to be
%
%
\begin{eqnarray} \label{eqproofXX2}
&& \frac{1}{2}\bigl\langle X(t),X(t)\bigr\rangle_{\caH_2}
\nonumber
\\
&&\quad = \int_0^t\int_0^s
\bigl\langle g(t,u)\,\intd B(u),g(t,s)\,\intd B(s)\bigr\rangle_{\caH_2}
\nonumber
\\
&&\quad = \int_0^t\int_0^s
\bigl\langle g(s,u)\,\intd B(u),g(s,s)\,\intd B(s)\bigr\rangle_{\caH_2}
\nonumber
\\
&&\qquad {}+ \int_0^t\int
_s^t\biggl\langle\int_0^r
g(r,u)\,\intd B(u),g(\intd r,s)\,\intd B(s)\biggr\rangle_{\caH_2}
\\
&&\qquad {} + \int_0^t\int
_0^s \biggl\langle\int_s^t
g(\intd v,u)\,\intd B(u),g(s,s)\,\intd B(s)\biggr\rangle_{\caH_2}
\nonumber
\\
&&\qquad {} - \int_0^t\int
_s^t \biggl\langle\int_s^r
g(t,u)\,\intd B(u),g(\intd r,s)\,\intd B(s)\biggr\rangle_{\caH_2}
\nonumber
\\
&&\qquad {} + \int_0^t\int
_s^t \biggl\langle\int_0^r
\bigl(g(t,u)-g(r,u)\bigr)\,\intd B(u),g(\intd r,s)\,\intd B(s)\biggr
\rangle_{\caH
_2}\nonumber
\end{eqnarray}
using the order $u\leq s\leq r\leq t$. Note that the first and second
term on the right-hand side of this equality are equal to first and
second term on the right-hand side of (\ref{eqproofXX}). The three
last terms in (\ref{eqproofXX2}) can be rewritten by grouping the
first and third one together as
\begin{eqnarray*}
&& \int_0^t\int_0^s
\biggl( \biggl\langle\int_s^t g(\intd v,u)\,\intd B(u), g(s,s)\,\intd B(s)\biggr\rangle_{\caH_2}
\\
&&\hspace*{25pt}\quad{}
+ \int_s^t \biggl\langle\int
_r^t g(\intd v,u)\,\intd B(u), g(\intd r,s)\,\intd
B(s)\biggr\rangle_{\caH_2} \biggr)
\\
&&\quad{} -\int_0^t\int_s^t
\biggl\langle\int_s^r g(r,u)\,\intd B(u),g(\intd
r,s)\,\intd B(s)\biggr\rangle_{\caH_2},
\end{eqnarray*}
which we can calculate to be, using Fubini's theorem,
\begin{eqnarray*}
&& \sum_{k,l,m\in\N} \int_0^t
\int_0^s \biggl(\biggl\langle\int
_s^tg(\intd v,u)e_l,e_k
\biggr\rangle_{\caH_2}\bigl\langle g(s,s)e_m,e_k
\bigr\rangle_{\caH_2}
\\
&&\hspace*{43pt}\qquad {} + \int_s^t\biggl\langle\int
_r^tg(\intd v,u)e_l,e_k
\biggr\rangle_{\caH_2}\bigl\langle g(\intd r,s)e_m,e_k
\bigr\rangle_{\caH_2} \biggr)\,\intd B^l(u)\,\intd
B^m(s)
\\
&&\qquad {} - \sum_{k,l,m\in\N} \int
_0^t\int_s^t
\int_s^r \bigl\langle g(r,u)e_l,e_k
\bigr\rangle_{\caH_2}\,\intd B^l(u)\bigl\langle g(\intd
r,s)e_m,e_k\bigr\rangle_{\caH_2}\,\intd
B^m(s)
\\
&&\quad = \sum_{k,l,m\in\N} \int_0^t
\int_0^s \caK_{\langle
g(e_l),e_k\rangle_{\caH_2}} \biggl(\int
_\cdot^t\bigl\langle g(\intd v,u)e_m,e_k
\bigr\rangle_{\caH_2} \biggr) (t,s)\,\intd B^l(u)\,\intd
B^m(s)
\\
&&\qquad {} - \sum_{k,l,m\in\N} \int
_0^t\int_s^t
\int_u^t \bigl\langle g(r,u)e_l,e_k
\bigr\rangle_{\caH_2}\bigl\langle g(\intd r,s)e_m,e_k
\bigr\rangle_{\caH_2}\,\intd B^l(u)\,\intd B^m(s)
\\
&&\quad = \int_0^t\int_0^s
\caK_g \biggl(\int_\cdot^t g(\intd
v,u)^* \biggr) (t,s)\,\intd B(u)\otimes\intd B(s)
\\
&&\qquad {} - \int_0^t\int
_s^t\int_u^t
\bigl\langle g(r,u),g(\intd r,s)\bigr\rangle_{\caH_2}\,\intd B(u)\otimes
\intd B(s),
\end{eqnarray*}
where $g(\cdot,s)^*$ means the adjoint linear operator of $g(\cdot
,s)$. The third term in (\ref{eqproofXX}) can be seen to be equal to
%
%
\begin{eqnarray}
&& \sum_{k,l\in\N}\int_0^t
\biggl(\bigl\langle g(s,s)e_l,e_k\bigr\rangle
_{\caH_2}\bigl\langle g(s,s)e_l,e_k\bigr
\rangle_{\caH_2} + \int_s^t \bigl\langle
g(r,s)e_l,e_k\bigr\rangle_{\caH_2}\bigl\langle
g(\intd r,s)e_l,e_k\bigr\rangle_{\caH_2}
\biggr)\,\intd s
\nonumber
\\
&&\quad = \sum_{k,l\in\N}\int_0^t
\caK_{\langle g(e_l),e_k\rangle_{\caH
_2}} \bigl(\bigl\langle g(\cdot,s)e_l,e_k
\bigr\rangle_{\caH_2} \bigr) (t,s)\,\intd s \label{eqintegrationXX}
\\
&&\quad = \tr_{\caH_2}\int_0^t
\caK_g\bigl(g(\cdot,s)^*\bigr) (t,s)\,\intd s.
\nonumber
\end{eqnarray}
One could also reformulate the integral term in the first line of
(\ref{eqintegrationXX}), which has the form $\int_s^t f(u)f(\intd u)
= \frac{1}{2}(f^2(t)-f^2(s))$, to be
\begin{eqnarray*}
&& \sum_{k\in\N}\int_0^t
\biggl(\bigl\langle g(s,s)e_l,e_k\bigr
\rangle_{\caH
_2}\bigl\langle g(s,s)e_l,e_k\bigr
\rangle_{\caH_2} + \int_s^t \bigl\langle
g(r,s)e_l,e_k\bigr\rangle_{\caH_2}\bigl\langle
g(\intd r,s)e_l,e_k\bigr\rangle_{\caH_2}
\biggr)\,\intd s
\\
&&\quad = \frac{1}{2}\int_0^t \bigl(\bigl
\llVert g(t,s)\bigr\rrVert^2_{L_2(\caH_1)} + \bigl\llVert g(s,s)
\bigr\rrVert^2_{L_2(\caH_1)} \bigr)\,\intd s.
\end{eqnarray*}
So, collecting all the terms above we finally have an equality
\begin{eqnarray*}
&& \int_0^t X(s)\,\intd X(s)
\\
&&\quad = \frac{1}{2}\bigl\langle X(t),X(t)\bigr\rangle_{\caH_2} - \int
_0^t\int_0^s
\caK_g \biggl(\int_\cdot^t g(\intd
v,u)^* \biggr) (t,s)\,\intd B(u)\otimes\intd B(s)
\\
&&\qquad{} - \int_0^t\int_s^t
\int_u^t \bigl\langle g(r,u),g(\intd r,s)\bigr
\rangle_{\caH_2}\,\intd B(u)\otimes\intd B(s) + \tr_{\caH_2}\int
_0^t \caK_g\bigl(g(\cdot,s)^*
\bigr) (t,s)\,\intd s,
\end{eqnarray*}
which gives us the desired link between the $X$-integral of $X$ and
$\llVert
X(t)\rrVert ^2_{\caH_2}$.
\end{pf}
%

\section{An It\^o formula}\label{secito}
In this section, we derive It\^o formulas for the processes $X$ and $Z$
defined in (\ref{eqdefX}) and (\ref{eqdefZ}). In order to do this,
we rely on the It\^o formula in Hilbert spaces with anticipating
integrands in \cite{grorudpardoux}, Proposition 4.10.

%
\begin{proposition}\label{propitogeneral}
Let $\caH,\caK$ be separable Hilbert spaces and $B$ be a $\caH
$-valued cylindrical Wiener process. Let moreover $F\in\caC^2(\caK
;\R)$ (the twice Fr\'echet differentiable functionals) and let
$(V(t))_{t\geq0}$ be the stochastic process defined by
\[
V(t) = V(0) + \int_0^t A(s)\,\intd s + \int
_0^t C(s)\delta B(s),
\]
where $V(0)\in\D^{1,4}(\caK)$, $A\in\L^{1,4}(\caK)$ and $C\in\L
^{2,p}(\caH,\caK)$ for some $p>4$, see (\ref{eqdefL}). Then
%
\begin{eqnarray*}
F\bigl(V(t)\bigr) &=& F\bigl(V(0)\bigr)  + \int_0^t
F'\bigl(V(s)\bigr)A(s)\,\intd s + \int_0^t
F'\bigl(V(s)\bigr)C(s)\delta B(s)
\\
&&{} + \frac{1}{2}\tr_{\caH}\int_0^t
F''\bigl(V(s)\bigr) \bigl(D^-V\bigr) (s) \bigl(C(s)
\bigr)\,\intd s
\\
&&{}+ \frac{1}{2}\tr_{\caH}\int_0^t
F''\bigl(V(s)\bigr) \bigl(C(s)\bigr) \bigl(C(s)\bigr)\,\intd s,
\end{eqnarray*}
where
%
\begin{eqnarray*}
\bigl(D^-V\bigr) (s) & =& 2D_sV(0) + 2\int_0^s
D_sA(r)\,\intd r + 2\int_0^s
D_sC(r)\delta B(r).
\end{eqnarray*}
\end{proposition}

We give some remarks on the definition of the various terms in the
above formula. Since $F'\in L(\caK,L(\caK,\R))$, the first two
integral terms take values in $\R$, as $F(V(0))$ does. For the fourth
term on the right-hand side, we use the definition of $F''$ to be an
element of $L(\caK,L(\caK,L(\caK,\R)))$, whereas for the fifth term
we use the equivalent formulation $F''\in L(\caK\otimes\caK,L_2(\caK
,\R))$. Then, since $(D^-V)(s)\in L_2(\caH,\caK)$ and $C(s)\in
L_2(\caH,\caK)$, we have that $F''(V(s))(D^-V)(s)(C(s))\in L_2(\caH
,L_2(\caH,\R))$ and $F''(V(s))(D^-V)(s)(C(s))\in L_2(\caH,L_2(\caH
,\R))$ so that the trace over $\caH$ is in both cases well-defined.
Note that this trace is our way of writing the inner products in the
original formula in \cite{grorudpardoux}, Proposition 4.10.

The first issue in this section is to extend the above proposition to
functions $F\in\caC^2(\caK;\caK_1)$ where $\caK_1$ is another
separable Hilbert space, not necessarily $\R$. We can reduce this to
applying Proposition \ref{propitogeneral} coordinatewise in the
following sense: set $F^k:=\langle F,e_k\rangle_{\caK_1}$ and observe
\begin{eqnarray*}
&& F^k\bigl(V(t)\bigr)
\\
&&\quad = F^k\bigl(V(0)\bigr) + \int_0^t
\bigl(F^k\bigr)'\bigl(V(s)\bigr)A(s)\,\intd s + \int
_0^t \bigl(F^k
\bigr)'\bigl(V(s)\bigr)C(s)\delta B(s)
\\
&&\qquad{} + \frac{1}{2}\tr_{\caH}\int_0^t
\bigl(F^k\bigr)''\bigl(V(s)\bigr)
\bigl(D^-V\bigr) (s)C(s)\,\intd s + \frac{1}{2}\tr_{\caH}\int
_0^t \bigl(F^k
\bigr)''\bigl(V(s)\bigr) \bigl(C(s)\bigr)C(s)\,\intd s
\\
&&\quad = \biggl\langle F
\bigl(V(0)\bigr) + \int_0^t F'
\bigl(V(s)\bigr)A(s)\,\intd s + \int_0^t
F'\bigl(V(s)\bigr)C(s)\delta B(s)
\\
&&\hspace*{3pt}\qquad{} + \frac{1}{2}\tr_{\caH}\int_0^t
F''\bigl(V(s)\bigr) \bigl(D^-V\bigr) (s)C(s)\,\intd s +
\frac{1}{2}\tr_{\caH}\int_0^t
F''\bigl(V(s)\bigr) \bigl(C(s)\bigr)C(s)\,\intd
s,e_k \biggr\rangle_{\caK_1}.
\end{eqnarray*}
This holds true since one can commute Fr\'echet derivatives and
projections, as they are bounded linear operators. Hence, we have
identified each coordinate and summing over these coordinates yields an
It\^o formula in the Hilbert-valued case.

Now we derive an It\^o formula for processes like $Z$ in (\ref
{eqdefZ}) which may be Hilbert-valued or real-valued. In order to apply
Proposition \ref{propitogeneral} to our case, we will have to assume
that $g(s,s)$ is a well-defined linear operator for all $s\in[0,T]$
and that there exists a function $\phi$ as in (\ref
{eqdifferentiableg}). Therefore, the It\^o formula in this section will
only hold for the case when $X$ is a Skorohod semimartingale, that is,
with $g$ satisfying the conditions of Proposition \ref
{propsmgcondition} but with an anticipating integrand in the
stochastic integral. So we will apply the above proposition to the
following stochastic process:
\begin{eqnarray*}
&& \int_0^tY(s)\,\intd X(s)
\\
&&\quad= \int_0^t \caK_g(Y) (t,s)
\sigma(s)\delta B(s) + \int_0^t \tr
_{\caH_1} \bigl(D_s\caK_g(Y) (t,s)\sigma(s)
\bigr)\,\intd s
\\
&&\quad = \int_0^t Y(s)g(s,s)\sigma(s)\delta B(s) +
\int_0^t\int_s^t
Y(u)\phi(u,s)\,\intd u\sigma(s)\delta B(s)
\\
&&\qquad {}+ \int_0^t
\tr_{\caH_1} D_s\bigl(Y(s)\bigr)g(s,s)\sigma(s)\,\intd s + \int
_0^t \tr_{\caH_1} \int
_s^t D_u\bigl(Y(u)\bigr)\phi(u,s)\,\intd u\sigma(s)\,\intd s
\\
&&\quad = \int_0^t b(s)\delta B(s)
+ \int_0^t \biggl(\int
_0^s a(s,u)\delta B(u) + \tr_{\caH_1}
\bigl(b_D(s) \bigr) + \tr_{\caH_1} \int_0^s
a_D(s,u)\,\intd u \biggr)\,\intd s,
\end{eqnarray*}
where we have used the stochastic Fubini theorem and
\begin{eqnarray*}
b(s) & =& Y(s)g(s,s)\sigma(s), \qquad a(s,u) = Y(s)\phi(s,u)\sigma(u),
\\
b_D(s) & =& D_s\bigl(Y(s)\bigr)g(s,s)\sigma(s), \qquad
a_D(s,u) = D_s\bigl(Y(s)\bigr)\phi(s,u)\sigma(u).
\end{eqnarray*}
So we apply Proposition \ref{propitogeneral} to the processes $C(s) =
b(s)$ and
%
%
\begin{equation}
A(s) = \int_0^s a(s,u)\delta B(u) +
\tr_{\caH_1} b_D(s) + \tr_{\caH
_1} \int
_0^s a_D(s,u)\,\intd u.
\label{eqAandC}
\end{equation}
After the following theorem which sets up an It\^o formula for the
integral with respect to $X$, we will provide some sufficient
conditions so that these processes satisfy the conditions in
Proposition \ref{propitogeneral}.

%
\begin{theorem}\label{thmito}
Let $\caH_4$ be a separable Hilbert space and let $F\dvtx \caH
_3\rightarrow\caH_4$ be twice Fr\'echet differentiable. Let $X$ and
$Z$ be defined as in (\ref{eqdefX}) and (\ref{eqdefZ}) where we
suppose that Assumption \ref{assdefintegral} holds. Furthermore,
assume that $g$ satisfies the semimartingale conditions in Proposition
\ref{propsmgcondition}. Assume that $Y$ and $\sigma$ are twice
Malliavin differentiable and
%
%
\begin{equation}
C(s)\in\L^{2,p}(\caH_1,\caH_3)\quad\mbox{and}
\quad A(s)\in\L^{1,4}(\caH_3), \label{eqconditionito}
\end{equation}
where $A(s)$ and $C(s)$ are as in (\ref{eqAandC}). Then $F'(Z)Y\in
\caI^X(0,t)$ for all $t\in[0,T]$ and
%
%
\begin{eqnarray}
\label{eqitoformula} F\bigl(Z(t)\bigr) &=& F(0)  + \int_0^t
\caK_g \bigl(F'(Z)Y \bigr) (t,s)\sigma(s)\delta B(s)
\nonumber
\\
&&{} + \tr_{\caH_1} \int_0^t
D_s \bigl(\caK_g \bigl(F'(Z)Y \bigr)
(t,s) \bigr)\sigma(s)\,\intd s
\nonumber
\\
&&{} - \frac{1}{2}\tr_{\caH_1}\int_0^t
F''\bigl(Z(s)\bigr) \bigl(Y(s)g(s,s)\sigma(s) \bigr)
\bigl(Y(s)g(s,s)\sigma(s) \bigr)\,\intd s
\\
&=& F(0)  + \int_0^t F'
\bigl(Z(s)\bigr)Y(s)\,\intd X(s)
\nonumber
\\
&&{} - \frac{1}{2}\tr_{\caH_1}\int_0^t
F''\bigl(Z(s)\bigr) \bigl(Y(s)g(s,s)\sigma(s) \bigr)
\bigl(Y(s)g(s,s)\sigma(s) \bigr)\,\intd s.\nonumber
\end{eqnarray}
\end{theorem}

\begin{pf}
Applying Proposition \ref{propitogeneral} yields
%
%
\begin{eqnarray} \label{eqproofito1}
F\bigl(Z(t)\bigr) &= & F(0) + \int_0^t
F'\bigl(Z(s)\bigr)b(s)\delta B(s) + \int_0^t
F'\bigl(Z(s)\bigr)\int_0^s
a(s,u)\delta B(u)\,\intd s
\nonumber
\\
&&{} + \int_0^t F'\bigl(Z(s)\bigr)
\tr_{\caH_1} b_D(s)\,\intd s + \int_0^t
F'\bigl(Z(s)\bigr)\tr_{\caH_1}\int_0^s
a_D(s,u)\,\intd u\,\intd s
\nonumber
\\
&&{} + \tr_{\caH_1} \int_0^t
F''\bigl(Z(s)\bigr)\int_0^s
D_sb(u)\delta B(u) b(s)\,\intd s
\\
&&{} + \tr_{\caH_1} \int_0^t
F''\bigl(Z(s)\bigr)\int_0^s
D_s\int_0^u a(u,r)\delta B(r)\,\intd u\, b(s)\,\intd s\nonumber
\\
&&{} + \tr_{\caH_1} \int_0^t
F''\bigl(Z(s)\bigr)\int_0^s
D_s \tr_{\caH_1} b_D(u)\,\intd u\,b(s)\,\intd s
\nonumber
\\
&&{} + \tr_{\caH_1} \int_0^t
F''\bigl(Z(s)\bigr)\int_0^s
D_s \tr_{\caH_1} \int_0^u
a_D(u,r)\,\intd r\,\intd u\,b(s)\,\intd s
\nonumber
\\
&&{} + \frac{1}{2}\tr_{\caH_1} \int_0^t
F''\bigl(Z(s)\bigr)b(s)b(s)\,\intd s.\nonumber
\end{eqnarray}
Note that the third term on the right-hand side of (\ref
{eqproofito1}) can be rewritten using (\ref{eqHIbP}) as
%
%
\begin{eqnarray}
\label{eqproofito2} && \int_0^t F'
\bigl(Z(s)\bigr)\int_0^s a(s,u)\delta B(u)\,\intd s
\nonumber\\[-8pt]\\[-8pt]\nonumber
&&\quad = \int_0^t \int_0^s
F'\bigl(Z(s)\bigr)a(s,u)\delta B(u)\,\intd s +\int
_0^t \int_0^s
\tr_{\caH_1} \bigl(\bigl(D_uF'\bigl(Z(s)\bigr)
\bigr) a(s,u) \bigr)\,\intd u\,\intd s.
\end{eqnarray}
By the It\^o formula in Proposition \ref{propitogeneral}, we know
that all the terms on the right-hand side of (\ref{eqproofito1}) are
well-defined. Now we calculate the terms on the right-hand side of
(\ref{eqitoformula}) \emph{assuming} that they exist. We will see
that they are equal to a sum of terms on the right-hand side of (\ref
{eqproofito1}), and we already know that the latter terms are
well-defined. Therefore, we can conclude that all the terms in (\ref
{eqitoformula}) are also well-defined and that the equality in (\ref
{eqitoformula}) holds. In fact, we have by definition
%
%
\begin{eqnarray}\label{eqproofito3}
&& \int_0^t \caK_g
\bigl(F'(Z)Y\bigr) (t,s)\sigma(s)\delta B(s)
\nonumber\\[-8pt]\\[-8pt]\nonumber
&&\quad = \int_0^t F'\bigl(Z(s)
\bigr)b(s)\delta B(s) + \int_0^t\int
_0^s F'\bigl(Z(s)\bigr)a(s,u)
\delta B(u)\,\intd s
\end{eqnarray}
and
%
%
\begin{eqnarray} \label{eqproofito4}
&& \tr_{\caH_1} \int_0^t
\bigl(D_s\caK_g\bigl(F'(Z)Y\bigr) (t,s)
\bigr)\sigma(s)\,\intd s
\nonumber
\\
&&\quad = \tr_{\caH_1} \int_0^t \biggl(
D_sF'\bigl(Z(s)\bigr)b(s)\,\intd s + F'
\bigl(Z(s)\bigr)b_D(s)\nonumber
\\
&&\hspace*{39pt}\qquad {} + \int_s^t
\bigl( \bigl(D_sF'\bigl(Z(u)\bigr) \bigr)a(u,s)+ F'\bigl(Z(u)\bigr)a_D(u,s) \bigr)\,\intd u \biggr)\,\intd s
\\
&&\quad = \tr_{\caH_1} \int_0^t
\bigl(D_sF'\bigl(Z(s)\bigr) \bigr)b(s)\,\intd s + \tr
_{\caH_1} \int_0^t F'
\bigl(Z(s)\bigr)b_D(s)\,\intd s
\nonumber
\\
&&\qquad {}+ \tr_{\caH_1}\int_0^t
\int_0^s \bigl(D_uF'
\bigl(Z(s)\bigr) \bigr)a(s,u)\,\intd u\,\intd s
+ \tr_{\caH_1}\int_0^t
\int_0^s F'\bigl(Z(s)
\bigr)a_D(s,u)\,\intd u\,\intd s.\nonumber
\end{eqnarray}
Note that by the chain rule (\ref{eqchainrule})
\begin{eqnarray*}
D_sF'\bigl(Z(s)\bigr)
&=& F''\bigl(Z(s)\bigr) \biggl(b(s) + \int
_0^s D_sb(u)\delta B(u) +
\tr_{\caH
_1} \int_0^s D_s
\int_0^u a(u,r)\delta B(r)\,\intd u
\\
&&\hspace*{34pt}\qquad {} + \int_0^sD_s
\tr_{\caH_1} b_D(u)\,\intd u + \int_0^s
D_s\tr_{\caH_1} \int_0^u
a_D(u,r)\,\intd r\,\intd u \biggr).
\end{eqnarray*}
With this, we can rewrite the first term on the right-hand side of
(\ref{eqproofito4}) to
%
%
\begin{eqnarray}
&& \tr_{\caH_1} \int_0^t
\bigl(D_sF'\bigl(Z(s)\bigr) \bigr)b(s)\,\intd s
\nonumber
\\
&&\quad = \tr_{\caH_1} \int_0^t
F''\bigl(Z(s)\bigr) b(s)b(s)\,\intd s +
\tr_{\caH_1} \int_0^t
F''\bigl(Z(s)\bigr)\int_0^sD_sb(u)
\delta B(u)b(s)\,\intd s
\nonumber
\\
&&\qquad {} + \tr_{\caH_1} \int_0^t
F''\bigl(Z(s)\bigr)\int_0^sD_s
\int_0^ua(u,r)\delta B(r)\,\intd u\,b(s)\,\intd s
\label{eqproofito5}
\\
&&\qquad {} + \tr_{\caH_1} \int_0^t
F''\bigl(Z(s)\bigr)\int_0^sD_s
\tr_{\caH_1} b_D(u)\,\intd u\,b(s)\,\intd s
\nonumber
\\
&&\qquad {} + \tr_{\caH_1} \int_0^t
F''\bigl(Z(s)\bigr)\int_0^sD_s
\tr_{\caH_1}\int_0^u
a_D(u,r)\,\intd r\,\intd u\,b(s)\,\intd s.
\nonumber
\end{eqnarray}
Now we collect terms. The two terms on the right-hand side of (\ref
{eqproofito3}) are the second term in (\ref{eqproofito1}) and the
first term in (\ref{eqproofito2}). Moreover, the last three terms on
the right-hand side of (\ref{eqproofito4}) are equal to the fourth
term in (\ref{eqproofito1}), the second term in (\ref{eqproofito2})
and the fifth term in (\ref{eqproofito1}). Furthermore, the last four
terms in (\ref{eqproofito5}) are the sixth to ninth term in (\ref
{eqproofito1}). Finally, the first term one the right-hand side of
(\ref{eqproofito5}) is twice the last term in (\ref{eqproofito1}).
Hence, by adding and subtracting the last term in (\ref{eqproofito1})
to this equality, we obtain the assertion.
\end{pf}

Next, we provide a sufficient condition under which (\ref
{eqconditionito}) holds relying on H\"older's inequality.

%
\begin{remark}
For the term $C(s)$ in (\ref{eqconditionito}), we have to deal with
the following three terms which constitute its $\L^{2,p}$-norm. We
will apply H\"older's Inequality with respect to $\omega$ (note that
$g$ is assumed to be deterministic) with the two conjugate exponents
$q_1,q_2\in[1,\infty]$ (where if one of them is infinite, the
supremum norm has to be used). We will now assume that $Y(s)\in
L_2(\caH_2,\caH_3)$, but a similar calculation would also be possible
with $\sigma(s)\in L_2(\caH_1,\caH_2)$. We obtain
\begin{eqnarray*}
&& \int_0^T \E\bigl[\bigl\llVert Y(s)g(s,s)
\sigma(s)\bigr\rrVert_{L_2(\caH_1,\caH
_3)}^p \bigr] \,\intd s
\\
&&\quad \leq\int_0^T \bigl\llVert g(s,s)\bigr
\rrVert^p_{L(\caH_2,\caH_2)} \bigl(\E\bigl[\bigl\llVert Y(s)\bigr
\rrVert^{pq_1}_{L_2(\caH_2,\caH_3)} \bigr] \bigr)^{1/q_1} \bigl(\E
\bigl[\bigl\llVert\sigma(s)\bigr\rrVert_{L(\caH_1,\caH_2)}^{pq_2}
\bigr]
\bigr)^{1/q_2} \,\intd s,
\end{eqnarray*}
similarly for the Malliavin derivative using (\ref{eqproductrule})
\begin{eqnarray*}
&& \int_0^T\int_0^T
\E\bigl[\bigl\llVert D_t\bigl(Y(s)g(s,s)\sigma(s)\bigr)\bigr\rrVert
_{L_2(L_2(\caH_0,\caH_1),\caH_3)}^p \bigr] \,\intd s\,\intd t
\\
&&\quad \leq\int_0^T \bigl\llVert g(s,s)\bigr
\rrVert^p_{L(\caH_2,\caH_2)}
\\
&&\hspace*{15pt}\qquad {}\times\biggl( \bigl(\E\bigl[\bigl\llVert\sigma(s)\bigr
\rrVert_{L(\caH
_1,\caH_2)}^{pq_2} \bigr] \bigr)^{1/q_2}\int
_0^T \bigl(\E\bigl[\bigl\llVert
D_tY(s)\bigr\rrVert^{pq_1}_{L_2(\caH_1\otimes\caH_2,\caH_3)} \bigr]
\bigr)^{1/q_1}\,\intd t
\\
&&\hspace*{33pt}\qquad {} + \bigl(\E\bigl[\bigl\llVert Y(s)\bigr\rrVert
_{L_2(\caH_2,\caH
_3)}^{pq_1} \bigr] \bigr)^{1/q_2}\int
_0^T \bigl(\E\bigl[\bigl\llVert
D_t\sigma(s)\bigr\rrVert^{pq_2}_{L_2(\caH_1\otimes\caH_1,\caH_2)} \bigr]
\bigr)^{1/q_1}\,\intd t \biggr)\,\intd s,
\end{eqnarray*}
and finally for the second Malliavin derivative
\begin{eqnarray*}
&& \int_0^T\int_0^T
\int_0^T \E\bigl[\bigl\llVert
D_rD_t\bigl(Y(s)g(s,s)\sigma(s)\bigr)\bigr\rrVert
_{L_2(L_2(\caH_0,\caH_1)\otimes\caH_1,\caH_3)}^p \bigr] \,\intd s\,\intd
t\,\intd r
\\
&&\quad \leq\int_0^T \bigl\llVert g(s,s)\bigr
\rrVert^p_{L(\caH_2,\caH_2)}
\\
&&\hspace*{15pt}\qquad{} \times\biggl(\bigl(\E\bigl[\bigl\llVert\sigma(s)\bigr\rrVert_{L(\caH
_1,\caH_2)}^{pq_2}
\bigr] \bigr)^{1/q_2}\int_0^T\int
_0^T \bigl(\E\bigl[\bigl\llVert
D_rD_tY(s)\bigr\rrVert^{pq_1}_{L_2(\caH_1\otimes\caH_1\otimes\caH
_2,\caH_3)}
\bigr] \bigr)^{1/q_1}\,\intd t\,\intd r
\\
&&\hspace*{33pt}\qquad{} + 2\int_0^T \bigl(\E\bigl[\bigl\llVert
D_rY(s)\bigr\rrVert_{L_2(\caH_2,\caH
_3)}^{pq_1} \bigr]
\bigr)^{1/q_1}\,\intd r
\\
&&\hspace*{61pt}{}\times \int_0^T \bigl(\E
\bigl[\bigl\llVert D_t\sigma(s)\bigr\rrVert^{pq_2}_{L_2(L_2(\caH_0,\caH
_1)\otimes\caH_1,\caH
_2)}
\bigr] \bigr)^{1/q_2}\,\intd t
\\
&&\hspace*{32pt}\qquad{} + \bigl(\E\bigl[\bigl\llVert Y(s)\bigr\rrVert_{L_2(\caH_2,\caH_3)}^{pq_1}
\bigr] \bigr)^{1/q_1}
\\
&&\hspace*{61pt}{}\times \int_0^T\int
_0^T \bigl(\E\bigl[\bigl\llVert
D_rD_t\sigma(s)\bigr\rrVert^{pq_2}_{L_2(L_2(\caH_0,\caH_1)\otimes\caH
_1,\caH_2)}
\bigr] \bigr)^{1/q_2}\,\intd t\,\intd r \biggr)\,\intd s.
\end{eqnarray*}
From the last equality, we see that $Y$, $DY$, $DDY$, $\sigma$,
$D\sigma$ and $DD\sigma$ appear and have to satisfy some
integrability conditions with respect to the temporal and spatial
argument. Therefore, we conclude that a sufficient condition for the
conditions in (\ref{eqconditionito}) is that $Y$ and $\sigma$ belong
to a Sobolev-type space with respect to all arguments (including the
random argument), that is, in some $\L^{2,q}$-spaces for $q$
sufficiently large, possibly with a different $q$ for each of the two.
Moreover, $g(s,s)$ has to be in some $L^p$-space for $p>4$ with respect
to $\lambda| _{[0,T]}$.

In some special situations, we can further manipulate the terms in the
equality above. For instance, one can assume that $\sigma$ is
independent of the noise (what would imply its Malliavin derivatives to
be zero), or $\sigma$ might be independent of $Y$ in which case the
expectations involving $Y$ and $\sigma$ would factorize. An
investigation of the term $A(s)$ in (\ref{eqconditionito}) does not
yield to any further conditions, except that also the function
$(t,s)\mapsto\llVert \phi(t,s)\rrVert ^4_{L(\caH_2,\caH_2)}$ has to be
integrable with respect to $\lambda^2| _{[0,T]\times[0,T]}$.
\end{remark}

A first immediate application of Theorem \ref{thmito} is to derive an
It\^o formula for $X$. Note that due to (\ref
{eqresultdifference15}), $Y(s)\equiv\id_{\caH_2}$ and, therefore,
\begin{eqnarray*}
X(t) & =& \int_0^t\intd X(s)
 = \int_0^t \caK_g(
\id_{\caH_2}) (t,s)\sigma(s)\delta B(s) + \tr_{\caH_1} \int
_0^t \bigl(D_s\caK_g(
\id_{\caH_2}) (t,s) \bigr)\sigma(s)\,\intd s.
\end{eqnarray*}
Using Theorem \ref{thmito}, we get the following result.

%
\begin{corollary}\label{coritoX}
Suppose that $X$ is defined as in (\ref{eqdefX}), where the
conditions in Assumption \ref{assdefintegral} hold, and that $F\dvtx \caH
_2\rightarrow\caH_3$ is twice Fr\'echet differentiable. Assume
moreover that $\sigma$ is twice Malliavin differentiable and that $g$
has the form (\ref{eqdifferentiableg}). Then $F'(X)\in\caI^X(0,t)$
for all $t\in[0,T]$ and
\begin{eqnarray*}
F\bigl(X(t)\bigr) &=& F(0)  + \int_0^t
F'\bigl(X(s)\bigr)\,\intd X(s)
 - \frac{1}{2}\tr_{\caH_1}\int_0^t
F''\bigl(X(s)\bigr) \bigl(g(s,s)\sigma(s) \bigr)
\bigl(g(s,s)\sigma(s) \bigr)\,\intd s.
\end{eqnarray*}
\end{corollary}

Note that this is consistent with the It\^o formulas in the
one-dimensional case treated in \cite{alosnualart}, Theorems 1, 2. In
fact, in that setting we have $\sigma\equiv1$, $\caH_1=\caH_2=\R$,
$X(t) = \int_0^t g(t,s)\,\intd W(s)$ and $R(s):=\int_0^s
g(s,r)g(s,r)\,\intd r$. Therefore, Corollary \ref{coritoX} implies
%
%
\begin{eqnarray}
F\bigl(X(t)\bigr) &=& F(0)  + \int_0^t
\biggl(F'\bigl(X(s)\bigr)g(s,s) + \int_s^t
F'\bigl(X(u)\bigr)\phi(u,s)\,\intd u \biggr)\delta B(s)
\nonumber
\\
&&{} + \int_0^t \biggl(D_sF'
\bigl(X(s)\bigr)g(s,s) + \int_s^t
D_sF'\bigl(X(u)\bigr)\phi(u,s)\,\intd u \biggr)\,\intd s
\label{eqitoX}
\\
&&{} - \frac{1}{2}\int_0^t
F''\bigl(X(s)\bigr)g(s,s)g(s,s)\,\intd s.
\nonumber
\end{eqnarray}
Note that $D_sF'(X(u)) = F''(X(s))g(u,s)$ for $u\geq s$. Note
furthermore that
\[
\frac{\intd R(s)}{\intd s} = g(s,s)^2 + 2\int_0^s
g(s,u)\phi(s,u)\,\intd u.
\]
This and Fubini's theorem yield that the last three terms in (\ref
{eqitoX}) are equal to
\[
\frac{1}{2}\int_0^t
F''\bigl(X(s)\bigr) \biggl(g(s,s)^2 + 2
\int_0^s g(s,u)\phi(s,u)\,\intd u \biggr)\,\intd
s = \frac{1}{2}\int_0^tF''
\bigl(X(s)\bigr)\,\intd R(s).
\]
This yields
\[
F\bigl(X(t)\bigr) = F(0) + \int_0^t
F'\bigl(X(s)\bigr)\,\intd X(s) + \int_0^t
F''\bigl(X(s)\bigr)\,\intd R(s),
\]
which is the formula in \cite{alosnualart}, Theorems 1, 2.

A second application of Theorem \ref{thmito} is to calculate $Z^2$ in
the case $\caH_3=\R$. For this, we suppose the same conditions as in
Theorem \ref{thmito}, assume $\caH_3=\R$ and $F(x)=x^2$. Then
applying Theorem \ref{thmito} yields
%
%
\begin{eqnarray}\label{eqproofZ21}
\frac{1}{2}\bigl(Z(t)\bigr)^2
& =& \int_0^t Z(s)Y(s)g(s,s)\sigma(s)\delta
B(s) + \int_0^t\int_0^s
Z(s)Y(s)\phi(s,u)\sigma(u)\delta B(u)\,\intd s
\nonumber
\\[-2pt]
&& {} + \tr_{\caH_1} \int_0^t
\bigl(D_sZ(s) \bigr)Y(s)g(s,s)\sigma(s)\,\intd s
\nonumber
\\[-2pt]
&&{} + \tr_{\caH_1}\int_0^t
Z(s) \bigl(D_sY(s) \bigr)g(s,s)\sigma(s)\,\intd s
\nonumber\\[-8pt]\\[-8pt]\nonumber
&& {} + \tr_{\caH_1} \int_0^t
\int_0^s \bigl(D_sZ(s) \bigr)Y(s)
\phi(s,u)\sigma(u)\,\intd u\,\intd s
\nonumber
\\[-2pt]
&& {} + \tr_{\caH_1}\int_0^t
\int_0^s Z(s) \bigl(D_sY(s)
\bigr)\phi(s,u)\sigma(u)\,\intd u\,\intd s
\nonumber
\\[-2pt]
&& {} - \frac{1}{2}\int_0^t
\bigl\llVert Y(s)g(s,s)\sigma(s)\bigr\rrVert^2_{L_2(\caH_1,\caH
_2)}\,\intd s.\nonumber
\end{eqnarray}
If we now introduce the symbolic notation
\begin{eqnarray*}
\intd Z(s) &= & \caK_g(Y) (\cdot,s)\sigma(s)\delta B(s) +
\tr_{\caH
_1} D_s\caK_g(Y) (\cdot,s)\sigma(s)\,\intd
s
\\[-2pt]
&= & Y(s)g(s,s)\sigma(s)\delta B(s) + \int_0^s
Y(s)\phi(s,u)\sigma(u)\delta B(u)\,\intd s
\\[-2pt]
&& {} + \tr_{\caH_1} \bigl(D_sY(s) \bigr)g(s,s)
\sigma(s)\,\intd s + \tr_{\caH_1} \int_0^s
\bigl(D_sY(s) \bigr)\phi(s,u)\sigma(u)\,\intd u\,\intd s,
\end{eqnarray*}
then we see that the term $\int_0^t Z(s)\,\intd Z(s)$ accounts for the
first, second, fourth, fifth and sixth term (use the rule in
Proposition \ref{propHIbP}) in (\ref{eqproofZ21}). The third and
last term remain as correction terms and one has the formal result
%
%
\begin{eqnarray}\label{eqitoZ2}
\frac{1}{2}\bigl(Z(t)\bigr)^2 &= & \int_0^t
Z(s)\,\intd Z(s) + \tr_{\caH_1} \int_0^t
\bigl(D_sZ(s) \bigr)Y(s)g(s,s)\sigma(s)\,\intd s
\nonumber\\[-8pt]\\[-8pt]\nonumber
&&{} - \frac{1}{2}\int_0^t \bigl\llVert
Y(s)g(s,s)\sigma(s)\bigr\rrVert^2_{L_2(\caH_1,\caH
_2)}\,\intd s.
\end{eqnarray}
One could also alter the correction terms. For instance, one could also
expand the second term on the right-hand side of (\ref{eqitoZ2}) by
calculating $D_sZ(s)$. This would reverse the sign of the third term on
the right-hand side of (\ref{eqitoZ2}), but one would get other
correction terms.

\section{A random-field approach to the $X$-integral}\label
{secrandomfieldintegral}
In this section, we investigate whether the random-field approach
pioneered in \cite{walsh} gives a reasonable interpretation of the
$X$-integral, that is, whether we can derive an integral which has the
form $\int_0^t\int_{\mathbb{R}^d}Y(s,y)X(\intd s,\intd y)$, where
%
%
\begin{equation}
X(t,x) = \int_0^t\int_{\mathbb{R}^d}g(t,s;x,y)
\sigma(s,y)M(\delta s,\intd y). \label{eqdefXrf}
\end{equation}
Here, similar to the terms in (\ref{eqdefX}) and (\ref{eqdefZ}),
$g$ is a deterministic function, $\sigma$ and $Y$ are random fields
and $M$ is a martingale measure and the integral in (\ref{eqdefXrf})
is understood in the Walsh sense. In the first part of this section, we
quickly review the concept of Walsh integration and summarize the main
ideas of a minor generalization for anticipating integrands. In the
second part, we present the random-field $X$-integral in the special
case of homogeneous noise and show that this integral and the one
derived in Section~\ref{secintegral} coincide.

\subsection{Walsh integration}\label{secwalsh}
Let $(M_t(A);t\in[0,T],A\in\caB_b({\mathbb{R}^d}))$ be a worthy Gaussian
martingale measure, where $\caB_b$ are the bounded Borel sets. This
means that each $M_t(A)$ has a Gaussian distribution on $\R$. Fix the
filtration $(\caF_t)_{t\in[0,T]}$ as the one generated by the
martingale measure and augmented by the sets of probability zero. The
martingale measure takes values in $L^2(\Omega)$. Since the martingale
measure is a martingale in $t$ for all fixed sets $A\in\caB
_b({\mathbb{R}^d})$,
we can associate a quadratic covariance functional to it, denoted by
$Q_M(t,A,B)=\langle M(A),M(B)\rangle_t$. We can define a stochastic
integral of elementary processes $f(s,x,\omega)=1_{ (a,b] }(s)1_{ A
}(x)Z(\omega)$, where $0\leq a<b\leq T$, $A\in\caB
_b({\mathbb{R}^d})$ and $Z$ is a bounded $\caF_a$-measurable random variable,
with respect to the martingale measure by
\[
\int_0^t\int_{\mathbb{R}^d}f(s,y)M(
\intd s,\intd y) = \bigl(M_{t\wedge
b}(A)-M_{t\wedge a}(A) \bigr)Z,
\]
and for simple processes (linear combinations of elementary processes)
by an obvious linear combination. Then, as in It\^o integration, this
stochastic integral is extended to a larger class of integrands by
using the isometry
%
%
\begin{eqnarray}\label{eqwalshisometry}
&& \E\biggl[ \biggl(  \int_0^t\int
_{\mathbb{R}^d}f(s,y)M(\intd s,\intd y) \biggr)^2 \biggr]
\nonumber\\[-8pt]\\[-8pt]\nonumber
&&\quad  = \E\biggl[\int_0^t\int
_{\mathbb{R}^d}\int_{\mathbb
{R}^d}f(s,x)f(s,y)Q_M(
\intd s,\intd x,\intd y) \biggr] =: \llVert f\rrVert^2_0,
\end{eqnarray}
where we extended the covariance functional to a measure on
$[0,T]\times{\mathbb{R}^d}\times{\mathbb{R}^d}$. In order to be
able to do that, one
needs the fact that the martingale measure is assumed to be worthy.
This means that there exists a dominating measure $K_M$, that is a
positive definite measure with some regularity conditions for which
$\llvert Q_M([0,t],A,B)\rrvert \leq K_M([0,t],A,B)$ for all $t\in
[0,T]$ and $A,B\in
\caB_b({\mathbb{R}^d})$. In this case one gets an upper bound for the isometry
(\ref{eqwalshisometry})
\begin{eqnarray*}
&& \E\biggl[  \int_0^t \int
_{\mathbb{R}^d}\int_{\mathbb
{R}^d}f(s,x)f(s,y)Q_M(
\intd s,\intd x,\intd y) \biggr]
\\
&&\quad  \leq\E\biggl[\int_0^t\int
_{\mathbb{R}^d}\int_{\mathbb
{R}^d}\bigl\llvert f(s,x)\bigr
\rrvert\bigl\llvert f(s,y)\bigr\rrvert K_M(\intd s,\intd x,\intd y)
\biggr] =: \llVert f\rrVert^2_+.
\end{eqnarray*}
The norm $\llVert \cdot\rrVert _0$ is actually induced by an inner
product and one
defines two possible domains $\caP_+$ and $\caP_0$ of the stochastic
integral. The former one is defined to be the set of all functions $f$
with $\llVert f\rrVert _+<\infty$ and the latter one is defined as
the completion
of the space of simple processes with respect to the inner product
$\langle\cdot,\cdot\rangle_0$. Note that by \cite{walsh},
Proposition~2.3, $\caP_+$ is complete and the simple processes form a
dense subset in this space. $\caP_0$ on the other hand is a Hilbert
space, but does not coincide with the set of all processes such that
$\llVert
f\rrVert _0<\infty$. Note that for particular choices of $M$, there
might be
even Schwartz distributions in $\caP_0$. A particular example for this
construction is the case of spatially homogeneous noise, where
%
%
\begin{eqnarray}\label{eqhomogeneousnoise}
Q_M\bigl([0,t],A,B\bigr) & =& K_M\bigl([0,t],A,B\bigr)
\nonumber\\[-8pt]\\[-8pt]\nonumber
& =& \int_0^T\int_{\mathbb{R}^d}
\int_{\mathbb{R}^d}1_{ [0,t]
}(s)1_{ A }(y)1_{ B }(y-x)\,\intd y\,\Gamma(\intd x)\,\intd s,
\end{eqnarray}
for a nonnegative, nonnegative definite tempered measure $\Gamma$.
Then one can find a spectral representation of the $\llVert \cdot
\rrVert _0$-norm
so that in this case $\caP_0$ contains distributions. Note however,
that in this article we only use the integration concept by Walsh for
functions in $\caP_+$ since we are exclusively interested in the case
when $g$ is a function, as this is the case with ambit fields in
Example \ref{exambitfields}. However, in Section~\ref
{secSPDEconncetion} we will also admit distributions for $g$ which
calls for the Dalang integral introduced in \cite{dalang}.

For a correct treatment, we need to provide a certain extension of the
Walsh integral in the case of homogeneous noise so that it can handle
anticipating integrands. This is, however, not very difficult if one
keeps in mind the reformulation of the Walsh integral as a sum of
independent Brownian motions on a Hilbert space; see \cite{dalangquer}
or as the divergence operator of Malliavin calculus, see \cite{sanzsuess1}.

We define $\caH$ to be the completion of the simple functions on
${\mathbb{R}^d}
$ by the scalar product
\[
\langle f,g\rangle_\caH= \int_{\mathbb{R}^d}(f \ast g) (z)
\Gamma(\intd z) = \int_{\mathbb{R}^d}\int_{\mathbb{R}^d}f(y)g(y-z)\,\intd y\,\Gamma(\intd z),
\]
with $\Gamma$ as above and $\caH_T:=L^2([0,T];\caH)$. In contrast to
\cite{dalang} we do not need to introduce the spectral representation
of this scalar product since we need $f$ and $g$ to be functions. Then
the worthy martingale measure $M$ can be extended into an isonormal
Gaussian process $(M(h);h\in\caH_T)$ with respect to which we can do
Malliavin calculus. With this, we can define $\D^{1,2}(\caH)$ and the
important space $\L^{1,2}(\caH):=L^2([0,T];\D^{1,2}(\caH))$,
with norm
\[
\llVert f\rrVert^2_{\L^{1,2}(\caH)} = \E\bigl[\llVert f\rrVert
^2_{\caH_T} \bigr] + \E\bigl[\llVert Df\rrVert
^2_{\caH_T\otimes\caH_T} \bigr]<\infty.
\]
On this space the divergence operator $\delta\dvtx L^2(\Omega;\caH
_T)\rightarrow L^2(\Omega)$ is well-defined and continuous; see \cite
{nualart}, equation~(1.47). For simple processes $f$,
which are all in $\L
^{1,2}(\caH)$ the anticipating Walsh integral is given by
\begin{eqnarray*}
&& \int_0^t f(t,x)  M(\delta t,\intd x)
 = \sum_{j=1}^d Z_j
\bigl(M_{t\wedge b_j}(A_j)-M_{t\wedge
a_j}(A_j) \bigr)
+ \sum_{j=1}^d\langle DZ_j,1_{ [a_j,b_j] }1_{ A_j
}
\rangle_{\caH_T}.
\end{eqnarray*}
Then one can extend this integral to integrands in the space $\L
^{1,2}(\caH)$, where the Malliavin differentiable simple processes are
dense. This extension is done using the following isometry:
%
%
\begin{eqnarray}\label{eq2isometryanticipating}
\E\bigl[(u\cdot M)^2\bigr]
&=& \E\bigl[\llVert u\rrVert^2_{\caH_T} \bigr] + \E\bigl[
\llVert Du\rrVert^2_{\caH
_T\otimes\caH_T} \bigr]
\nonumber
\\
& =& \E\bigl[\llVert u\rrVert^2_{\caH_T} \bigr]
\nonumber\\[-8pt]\\[-8pt]\nonumber
&&{} + \E\biggl[
\int_0^t\int_0^t
\int_{\mathbb{R}^d}\int_{\mathbb{R}^d}\int
_{\mathbb{R}^d}\int_{\mathbb{R}^d}D_{s_2,y_2-z_2}u(s_1,y_1)\nonumber
\\
&&\hspace*{113pt} {} \times D_{s_1,y_1-z_1}u(s_2,y_2)\,\intd y_1\,\intd y_2\,\Gamma(\intd z_1)\Gamma(
\intd z_2)\,\intd s_1\,\intd s_2 \biggr].\nonumber
\end{eqnarray}
Another way to think about this integral is by using the equivalence of
the anticipating Walsh integral and an infinite sum of anticipating
integrals with respect to Brownian motion. In fact, we have for all
$u\in\L^{1,2}(\caH)$
\[
\int_0^t\int_{\mathbb{R}^d}u(s,y)M(
\delta s,\intd y) = \sum_{k\in
\N} \int
_0^t \bigl\langle u(s,\cdot),e_k
\bigr\rangle_\caH\delta B^k(s),
\]
where $(e_k)_{k\in\N}$ is a CONS of $\caH$, $B$ is a Brownian motion
on $\caH$ and $B^k:=\langle B,e_k\rangle_{\caH}$ are independent
real-valued Brownian motions, for which anticipating calculus is well
known. Therefore, one could also take the last equality as the
definition of the anticipating Walsh integral.

\subsection{The random-field $X$-integral}
Now we start with the definition of the $X$-integral using a
random-field approach. In contrast to the definition of the
$X$-integral in Section~\ref{secintegral}, we have to start here with
defining the stochastic integral on elementary processes first and then
extend it to simple processes and further. This is due to the fact that
a similar integration by parts procedure done in Section~\ref
{secintegral} is not easily applicable since one would get nontrivial
boundary terms which are hard to interpret.

Take $\caH$ and $\caH_T$ to be the Hilbert spaces defined in the
previous subsection. The martingale measure that we will use in the
following will be the one using the homogeneous noise in (\ref
{eqhomogeneousnoise}), which appeared in \cite{dalang}. We will
however not use it in its full generality, but we only look at
functions $f$ for which the following norm:
%
%
\begin{equation}
\llVert f\rrVert_0^2 = \int_0^T
\int_{\mathbb{R}^d}\int_{\mathbb
{R}^d}f(s,y)f(s,y-z)\,\intd y\,
\Gamma(\intd z)\,\intd s < \infty. \label{eqcovmeasure}
\end{equation}
This is a norm -- and in fact it is generated by an inner product
$\langle\cdot,\cdot\rangle_0$ -- as one can see by going to its
spectral representation. We will however not do this, since we cannot
deal with distributions as integrators $g$ in the kernel (\ref
{eqdefKg}). Therefore, we define the Hilbert space $\caH_T$ as the
completion of the simple functions with respect to the inner product
$\langle\cdot,\cdot\rangle_0$. This is then a function space which
does not contain distributions. We also look at random functions in
$L^2(\Omega;\caH_T)$ equipped with the obvious norm and inner
product. The following assumptions are made.

%
\begin{assumption}\label{assrandomfieldintegral}
Fix $T>0$ and $t\in[0,T]$, and let $M=(M_t(A);t\in[0,T],A\in\caB
_b({\mathbb{R}^d}))$ be a worthy martingale measure with covariation
measure as
in (\ref{eqcovmeasure}). Furthermore, $g\dvtx [0,T]^2\times({\mathbb{R}^d}
)^2\rightarrow\R$ is a deterministic function for all $0\leq s <
t\leq T$ and all $x,y\in{\mathbb{R}^d}$ and $(\sigma(t,x); (t,x)\in
[0,T]\times
{\mathbb{R}^d})$ is a Malliavin differentiable random field such that
$(s,y)\mapsto g(t,s,x,y)\*\sigma(s,y)$ is integrable with respect to $M$
for all $(t,x)\in[0,T]\times{\mathbb{R}^d}$, that is, this random
function is
in $\L^{1,2}(\caH)$. Assume for all $(s,y)\in[0,T]\times{\mathbb
{R}^d}$ that
the function $(t,x)\mapsto g(t,s,x,y)$ has bounded variation on
$[u,v]\times{\mathbb{R}^d}$ for all $0\leq s<u<v\leq t$.
\end{assumption}

Now we turn to the definition of the integral. To this end, let first
$0\leq a<b\leq T$, $A\in\caB_b({\mathbb{R}^d})$ be of the form
$A=\bigtimes
_{j=1}^d [a_j,b_j]$ and set $Y(t,x,\omega)=1_{ (a,b] }(t)1_{ A
}(x)Z(\omega)$, where $Z$ is a bounded Malliavin differentiable
random variable, not necessarily measurable with respect to $\caF_a$.
Note that since the spatial argument of the integrator process $X$ only
appears in $g$, and since $X$ is linear in the kernel $g$, one can
derive using (\ref{eqfunctionA})
%
%
\begin{eqnarray}\label{eqXt,A}
X(t,A) & =& \int_0^t\int_{\mathbb{R}^d}g(t,s,A,y)
\sigma(s,y)M(\intd s,\intd y)
\nonumber\\[-8pt]\\[-8pt]\nonumber
& =& \int_0^t\int_{\mathbb{R}^d}\int
_{\mathbb{R}^d}1_{ A }(z) g(t,s,\intd z,y)\sigma(s,y)M(\intd
s,\intd y).
\end{eqnarray}
In the calculation that follows, we aim at arriving at a similar kernel
as in (\ref{eqdefKg}). In order to keep the notation tidy, we do this
calculation for the case when $g(s,s;\cdot,\cdot)$ is a well-defined
object, thus yielding a kernel as in (\ref{eqdefKg2}). By following
the same arguments as given below, one can however also arrive at
(\ref{eqdefKg}). Using the obvious definition for the stochastic
integral for simple processes and (\ref{eqXt,A}), one derives for
the simple process $Y$
\begin{eqnarray*}
&& \int_0^t\int_{\mathbb{R}^d}Y(s,y)X(
\intd s,\intd y)
\\
&&\quad := Z \bigl(X(t\wedge b,A)-X(t\wedge a,A) \bigr)
\\
&&\quad  = Z\int_0^t\int_{\mathbb{R}^d}
\biggl(1_{ [a,b] }(s)\int_{\mathbb
{R}^d}1_{ A }(z)g(s,s;
\intd z,y)
\\
&&\hspace*{44pt}\qquad {} + 1_{ [0,b] }(s) \biggl(\int_{\mathbb{R}^d}1_{ A
}(z)g(b,s;
\intd z,y) - \int_{\mathbb{R}^d}1_{ A }(z)g(s,s;\intd z,y)
\biggr)
\\
&&\hspace*{44pt}\qquad {} - 1_{ [0,a] }(s) \biggl(\int_{\mathbb{R}^d}1_{ A
}(z)g(a,s;
\intd z,y)
\\
&&\hspace*{114pt}{}-\int_{\mathbb{R}^d}1_{ A }(z)g(s,s;\intd z,y)
\biggr) \biggr)\sigma(s,y)M(\delta s,\intd y)
\\
&&\quad = Z\int_0^t\int_{\mathbb{R}^d}
\biggl(1_{ [a,b] }(s)\int_{\mathbb
{R}^d}1_{ A }(z)g(s,s;
\intd z,y) + 1_{ [0,a] }(s)\int_{a}^{b}\int
_{\mathbb{R}^d}1_{ A }(z)g(\intd u,s;\intd z,y)
\\
&&\hspace*{44pt}\qquad {} + 1_{ [a,b] }(s)\int_{s}^{b}
\int_{\mathbb{R}^d}1_{ A
}(z)g(\intd u,s;\intd z,y) \biggr)
\sigma(s,y)M(\delta s,\intd y).
\end{eqnarray*}
Note that the last two integral terms one the right-hand side of the
previous equation are equal to
\[
\int_{a\vee s}^{b\wedge t}\int_{\mathbb{R}^d}1_{ A }(z)g(
\intd u,s;\intd z,y) = \int_s^t\int
_{\mathbb{R}^d}1_{ [a,b] }(u)1_{ A }(z)g(\intd u,s;\intd
z,y).
\]
Plugging this into the equation above yields
%
%
\begin{eqnarray} \label{eqderivationrandomfieldintegral}
&& \int_0^t\int_{\mathbb{R}^d}Y(s,y)X(
\intd s,\intd y)
\nonumber
\\
&&\quad = Z\int_0^t\int_{\mathbb{R}^d}
\biggl(\int_{\mathbb{R}^d}1_{
[a,b] }(s)1_{ A }(z)g(s,s;
\intd z,y)
\\
&&\hspace*{44pt}\qquad {} + \int_s^t
1_{ [a,b] }(u)1_{ A }(z)g(\intd u,s;\intd z,y) \biggr)
\sigma(s,y)M(\delta s,\intd y).\nonumber
\end{eqnarray}
Similar to (\ref{eqdefKg}), we define an integration kernel
%
%
\begin{eqnarray}\label{eqdefKgfield}
&& \caK_g(h) (t,s,y)
\nonumber\\[-8pt]\\[-8pt]\nonumber
&&\quad := \int_{\mathbb{R}^d}h(s,z)g(t,s;\intd z,y) + \int
_s^t\int_{\mathbb{R}^d}
\bigl(h(u,z)-h(s,z) \bigr)g(\intd u,s;\intd z,y).
\end{eqnarray}
Note that this term is well-defined if $g$ is a function of bounded
variation spatial argument $z$ and of bounded variation on subintervals
bounded away from $s$ in the temporal argument $t$. As in Section~\ref
{secintegral}, one can rewrite this kernel in some situations. If
$g(s,s;z,y)$ exists almost everywhere (with respect to $\intd y\,\Gamma
(\intd z)\,\intd s$), then one has
\[
\caK_g(h) (t,s,y) = \int_{\mathbb{R}^d}h(s,z)g(s,s;\intd
z,y) + \int_s^t\int_{\mathbb{R}^d}h(u,z)g(
\intd u,s;\intd z,y).
\]
If the function $g$ has a partial derivative with respect to each
coordinate, then
\begin{eqnarray*}
\caK_g(h) (t,s,y) &=& \int_{\mathbb{R}^d}h(s,z)
\frac{\partial
^d}{\partial
z_1\cdots\partial z_d}g(t,s;z,y)\,\intd z
\\
&&{} + \int_s^t\int_{\mathbb{R}^d}
\bigl(h(u,z)-h(s,z) \bigr)\frac
{\partial
^{d+1}}{\partial u\,\partial z_1\cdots\partial z_d}g(u,s;z,y)\,\intd z\,\intd u.
\end{eqnarray*}
Using the kernel in (\ref{eqdefKgfield}), we obtain by pulling the
random variable inside the stochastic integral
\begin{eqnarray*}
\int_0^t\int_{\mathbb{R}^d}Y(s,y)X(
\intd s,\intd y) & =& Z\int_0^t\int
_{\mathbb{R}^d}\caK_g(1_{ [a,b] }1_{ A
})
(t,s,y)\sigma(s,y)M(\delta s,\intd y)
\\
& =& \int_0^t\int_{\mathbb{R}^d}
\caK_g(1_{ [a,b] }1_{ A
}Z) (t,s,y)\sigma(s,y)M(
\delta s,\intd y)
\\
&& {} + \bigl\langle DZ,\caK_g(1_{ [a,b] }1_{ A
})
(t,\cdot,\ast)\sigma(\cdot,\ast)\bigr\rangle_{\caH_T}.
\end{eqnarray*}
The second term of the previous equality can be manipulated as follows
in order to get an expression which explicitly involves the original
integrand $Y$
%
%
\begin{eqnarray}\label{eqrepresentationpathwiseintegral}
&& \bigl\langle DZ,\caK_g(1_{ [a,b] }1_{ A }) (t,
\cdot,\ast)\sigma(\cdot,\ast)\bigr\rangle_{\caH_T}
\nonumber
\\
&&\quad= \int_0^t\int
_{\mathbb{R}^d}\int_{\mathbb
{R}^d}D_{s,y-z}Z
\caK_g(1_{ [a,b] }1_{ A }) (t,s,y)\sigma(s,y)\,\intd y\,
\Gamma(\intd z)\,\intd s
\\
&&\quad= \int_0^t\int
_{\mathbb{R}^d}\int_{\mathbb
{R}^d}D_{s,y-z}\caK
_g(Y) (t,s,y)\sigma(s,y)\,\intd y\,\Gamma(\intd z)\,\intd s.
\nonumber
\end{eqnarray}
In total, we have for the random-field $X$-integral
%
%
\begin{eqnarray}\label{eqdefintfield}
\int_0^t\int_{\mathbb{R}^d}Y(s,y)X(
\intd s,\intd y)  & =& \int_0^t
\int_{\mathbb{R}^d}\caK_g(Y) (t,s,y)\sigma(s,y)M(\delta s,
\intd y)
\nonumber\\[-8pt]\\[-8pt]\nonumber
&& {} + \int_0^t\int
_{\mathbb{R}^d}\int_{\mathbb
{R}^d}D_{s,y-z}\caK
_g(Y) (t,s,y)\sigma(s,y)\,\intd y\,\Gamma(\intd z)\,\intd s.
\nonumber
\end{eqnarray}
First, we see the similarities to the definition in (\ref
{eqdefintegral}), where here we add the spatial component which in turn
forces us to include a possible spatial correlation into the definition
of the integral. Note that if we consider the uncorrelated case, then
$\Gamma=\delta_0$ and one gets an expression without~$z$. Now we
extend the definition of the integral by using the isometry in (\ref
{eq2isometryanticipating}) so that the integral (\ref{eqdefintfield})
is defined for all $Y\in\L^{1,2}(\caH)$ for which similar
integrability conditions as in Definition~\ref{defkernel} hold.

%
\begin{definition}\label{defrandomfieldointegral}
Let $X$ be the random field defined in (\ref{eqdefXrf}) together with
Assumption~\ref{assrandomfieldintegral}. We say that the random field
$Y=(Y(t,x);(t,x)\in[0,T]\times{\mathbb{R}^d})\in\caP_0$ belongs to
the domain
of the stochastic integral with respect to the random field $X$, if:
\begin{longlist}[(iii)]
\item[(i)] the process $(Y(s,z))_{z\in{\mathbb{R}^d}}$ is
integrable
with respect to $g(t,s;\intd z,y)$ almost surely and $(t,s,y)$-almost
everywhere,

\item[(ii)] the process $(Y(u,z)-Y(s,z))_{u\in(s,t]\times
{\mathbb{R}^d}}$ is integrable with respect to $g(\intd u,s;\intd
z,y)$ almost
surely and $(s,y)$-almost everywhere,

\item[(iii)]  $(s,y)\mapsto\caK_g(Y)(t,s,y)\sigma
(s,y)1_{ [0,t] }(s)$ is in the domain of the martingale measure $M$,
that is, $\caK_g(Y)(t,\cdot,\ast)\sigma(\cdot,\ast)1_{ [0,t]
}(\cdot)\in\caP_0$ and

\item[(iv)] $\caK_g(Y)(t,s,y)$ is Malliavin
differentiable with respect to $D_{s,y-z}$ for all $s\in[0,t]$ and
$y,z\in{\mathbb{R}^d}$ and the random field
\[
(s,y,z)\mapsto\tr_{\caH_1} D_{s,y-z} \bigl(\caK_g(Y)
(t,s,y) \bigr)\sigma(s,y)
\]
$\lambda\llvert _{[0,T]}\otimes\lambda\rrvert _{\mathbb{R}^d}\otimes
\Gamma
$-integrable on
$[0,t]\times{\mathbb{R}^d}\times{\mathbb{R}^d}$ almost surely.
\end{longlist}
We denote this by $Y\in\caI^X([0,t]\times{\mathbb{R}^d})$ and the integral
$\int_0^t\int_{\mathbb{R}^d}Y(s,y)X(\intd s,\intd y)$ is defined by~(\ref{eqdefintfield}).
\end{definition}

The main reason why we have introduced this integral is to be able to
integrate processes $(Y(t,x,y);t\in[0,T],x,y\in{\mathbb{R}^d})$ with
respect to
$t$ and $y$ and obtain a random field which has a pointwise (in $x$)
interpretation as a real-valued random variable, rather than as an
element in some abstract Hilbert space. This could serve to deduce
properties of ambit fields such as path continuity, or existence of
densities at each point.

We could easily derive similar properties for the random-field
$X$-integral as in Section~\ref{secrules}, but we restrict ourselves
to showing that the random field integral can be rewritten as a
Hilbert-valued stochastic integral when we interpret the process $X$ as
a stochastic process with values in some Hilbert space, as it was shown
in \cite{dalangquer,sanzsuess1}.

%
\begin{proposition}
Let $Y\in\caI^X([0,t]\times{\mathbb{R}^d})$. Then $Y\in\caI
^X(0,t)$ with
$\caH_1=\caH$ and $\caH_3=\R$ and
\[
\int_0^t\int_{\mathbb{R}^d}Y(s,y)X(
\intd s,\intd y) = \int_0^t Y(s)\,\intd X(s),
\]
where $X$ is interpreted as a random field on the left-hand side and as
an $\caH_2$-valued Volterra process on the right-hand one.
\end{proposition}

\begin{pf}
At first we argue that if $Y\in\caI^X([0,t]\times{\mathbb{R}^d})$,
then $Y\in
\caI^X(0,t)$. In fact, from the remark at the end of the previous
subsection it follows that if $Y$ is in the domain of the martingale
measure integral, then it is also in the domain of the integral with
respect to a $\caH$-valued Brownian motion. Furthermore, the
properties Definition~\ref{defrandomfieldointegral}(ii)--(iv) are exactly Definition \ref{defkernel}(i)--(iii) written out in the special case of $\caH_1=\caH$.
This yields the first assertion.

Second, we show that the stochastic integrals with respect to the
martingale measure and the $\caH$-valued Brownian motion are
equivalent. This follows by a straightforward adaption of the proof of
\cite{dalangquer}, Proposition 2.6, because by Definition\ref
{defrandomfieldointegral}(iii) we conclude that
\[
\caK_g(Y) (t,\cdot,\ast)\sigma(\cdot,\ast)1_{ [0,t] }(\cdot)
\in L^2(\Omega;\caH_T)
\]
and
\begin{eqnarray*}
&& \int_0^t\int_{\mathbb{R}^d}
\caK_g(Y) (t,s,y)\sigma(s,y)M(\delta s,\intd y)
= \sum_{k\in\N}\int_0^t
\bigl\langle\caK_g(Y) (t,s,\ast)\sigma(s,\ast),e_k
\bigr\rangle_{\caH}\,\intd B^k(s),
\end{eqnarray*}
where $B$ is an $\caH$-valued Wiener process and $B^k:= \langle
B,e_k\rangle_\caH$. This yields the equality of the stochastic
integrals in the $X$-integral and the random-field $X$-integral.

For the pathwise integrals, we start with the elementary processes
given by $Y(t,x,\omega)=1_{ [a,b] }(t)1_{ A }(x)Z(\omega)$ and
use the representation in (\ref{eqrepresentationpathwiseintegral}) to obtain
\begin{eqnarray*}
&& \int_0^t\int_{\mathbb{R}^d}\int
_{\mathbb{R}^d}D_{s,y-z}\caK_g(Y) (t,s,y)\sigma
(s,y)\,\intd y\,\Gamma(\intd z)\,\intd s
\\
&&\quad  = \int_0^t\bigl\langle D_sZ,
\caK_g(1_{ [a,b] }1_{ A }) (t,s,\ast)\sigma(s,\ast)
\bigr\rangle_\caH\,\intd s
\\
&&\quad  = \sum_{k\in\N} \int_0^t
\langle D_sZ,e_k\rangle_\caH\bigl\langle
\caK_g(1_{ [a,b] }1_{ A }) (t,s,\ast)\sigma(s,
\ast),e_k\bigr\rangle_\caH\,\intd s
\\
&&\quad = \sum_{k\in\N} \int_0^t
\bigl\langle D_{s,k}\caK_g(Y) (t,s,\ast)\sigma(s,
\ast),e_k\bigr\rangle_\caH\,\intd s
\\
&&\quad = \tr_\caH\int_0^t
D_s\caK_g(Y) (t,s)\sigma(s)\,\intd s.
\end{eqnarray*}
Then we extend the last equality to all processes $Y\in\caP_0$ which
satisfy Definition \ref{defrandomfieldointegral}(iv).
This, together with the equality above, implies the assertion.
\end{pf}

Now we treat a first example for the random field $X$-integral, another
one will follow in Section~\ref{secSPDEconncetion}. Here, we focus on
the nonlinear stochastic heat equation with null initial condition
\begin{eqnarray*}
\biggl(\frac{\partial}{\partial t}(t,x) - \Delta\biggr)u(t,x) & =&
b\bigl(u(t,x)\bigr) +
\sigma\bigl(u(t,x)\bigr)\dot{F}(t,x)\qquad \bigl(t,x\in(0,T]\times{\mathbb
{R}^d}\bigr)
\\
u(0,x) & =& 0,\qquad x\in{\mathbb{R}^d}.
\end{eqnarray*}
The fundamental solution to the associated PDE (the heat equation) is
given in any spatial dimensions by
\[
g(t,x) = \frac{1}{(4\uppi t)^{d/2}}\exp\biggl(\frac{-\llvert x\rrvert
^2}{4t} \biggr).
\]
Let in the following without much loss of generality $d=1$. Note that
the random-field solution to the SPDE above is given by
\begin{eqnarray*}
u(t,x) &= & \int_0^t\int_{\mathbb{R}^d}g(t-s,x-y)
\sigma\bigl(u(s,y)\bigr)M(\intd s,\intd y)
\\
&&{} + \int_0^t\int_{\mathbb{R}^d}g(t-s,x-y)b
\bigl(u(s,y)\bigr)\,\intd y\,\intd s,
\end{eqnarray*}
almost surely for all $(t,x)\in[0,T]\times{\mathbb{R}^d}$, with $M$
being the
martingale measure corresponding to the noise $\dot{F}$. We see that
the stochastic integral has the form of a Volterra process and,
therefore, we can apply the integration theory developed in this
section. Since the fundamental solution $g$ is differentiable in both
arguments and stationary in time and space, we can calculate $\caK
_g(Y)(t,s,y)$ to be
\begin{eqnarray*}
&& \caK_g(Y) (t,s,y)
\\
&&\quad = \int_{\mathbb{R}^d}Y(s,y+z)g(t-s,\intd z) + \int
_0^{t-s}\int_{\mathbb{R}^d} Y(u+s,z+y)g(
\intd u,\intd z)
\\
&&\quad = \int_{\mathbb{R}^d}Y(s,y+z)\partial_z g(t-s,z)\,\intd z
+ \int_0^{t-s}\int_{\mathbb{R}^d}Y(u+s,z+y)
\partial^2_{(u,z)}g(u,z)\,\intd z\,\intd u,
\end{eqnarray*}
where $\partial_z$ denotes the partial derivative with respect to the
spatial argument, and $\partial^2_{(u,z)}$ is the mixed spatial and
temporal argument. Note the heat kernel estimates for the derivatives
\[
\bigl\llvert\partial_x\,\partial^a_t
g(t,s,x,y)\bigr\rrvert\leq c(x-y) (t-s)^{-(3+2a)/2}\exp\biggl(-c
\frac{(x-y)^2}{t-s} \biggr)
\]
for $a\in\{0,1\}$. This implies that all $Y$ which are almost surely
polynomially bounded are integrable with respect to $g$. This implies (i)
and (ii) in Definition \ref
{defrandomfieldointegral}. In order for (iii) to be
satisfied, one needs that $\caK_g(Y)(t,\cdot,\ast)\sigma(\cdot
,\ast)\in\caP_0$. This depends on the concrete form of $\sigma$ and~$\Gamma$. In order for the pathwise integral to well-defined in
$L^2(\Omega)$ which implies (iv), one has to assume that
$Y$ is Malliavin differentiable and that the integrand of the pathwise
integral is in $L^2([0,T]\times\Omega;L^1({\mathbb{R}^d}))$.

\section{Generalization of the random-field $X$-integral}\label{secSPDEconncetion}
In this section, we follow two objectives. First, we want to present
another explicit example for the $X$-integral defined in Section~\ref
{secintegral}, and thus show how nonregular fundamental solutions to
partial differential equations enter into this framework. Second, we
want to show a first idea how to generalize the random-field
$X$-integral in Section~\ref{secrandomfieldintegral} in a way similar
to the generalizations of the Walsh integral by Dalang in \cite
{dalang}, going from functions to distributions as integrands. More
specifically, we investigate under which conditions on $Y$ the kernel
$\caK_g(Y)(t,s)$ is well-defined if $g$ is chosen to be the
fundamental solution to a partial differential operator which is
singular in the spatial argument.

For this, we focus on the specific example of a wave equation in
different spatial dimensions $d$ given by
%
%
\begin{equation}
\biggl(\frac{\partial^2}{\partial t^2}-\Delta_d \biggr)u(t,x) =
\delta_{0,0}, \label{eqwaveequation}
\end{equation}
with null initial conditions, where $\Delta_d$ is the $d$-dimensional
Laplace operator and $\delta_{0,0}$ is the Dirac delta distribution in
time and space. The solution to this equation is the fundamental
solution to the wave equation, which differs with the spatial
dimension. In the dimensions $d=1$ and $d=2$ the fundamental solution
is a function, which we can treat with similar methods as in the
example with the heat kernel in Section~\ref{secrandomfieldintegral}.

For $d=3$, the fundamental solution is given by $g(t,s)=c\rho
^3_{t-s}/(t-s)$, where $c>0$ is a constant and $\rho_t^3$ is the
surface measure on the sphere in three dimensions with radius $t$.
Therefore, the mass of this measure is equal to $ct^2$. Since this
fundamental solution is not a regular function anymore, the
random-field $X$-integral of Section~\ref{secrandomfieldintegral} is
no longer well-defined. So we use the Hilbert-valued integral from
Section~\ref{secintegral} with the Hilbert space $\caH$ introduced
in Section~\ref{secwalsh}.

%
\begin{proposition}
Suppose that Assumption \ref{assdefintegral} and the integrability
conditions for $\sigma$ in Definition \ref{defkernel} hold. Assume
that $Y$ be a random element with values in the linear functionals on
$\caH$ such that $Y$ is differentiable and $Y'$ does not have a
singularity of order greater or equal than $2$ at zero, that is, $Y'\in
\mathrm{o}(v^{-2})$ at zero. Then, for the special choice of
$g(t,s)=c\rho^3_{t-s}/(t-s)$, $Y\in\caI^X(0,t)$.
\end{proposition}

\begin{pf}
By performing integration by parts, we see that
\begin{eqnarray*}
c^{-1}\caK_g(Y) (t,s) 
& =& Y(s)
\frac{\rho^3_{t-s}}{t-s} + \int_0^{t-s}
\bigl(Y(v+s)-Y(s) \bigr)\frac{\intd}{\intd v} \biggl(\frac{\rho
^3_v}{v} \biggr)\,\intd
v
\\
& =& Y(s)\frac{\rho^3_{t-s}}{t-s} + \bigl(Y(v+s)-Y(s) \bigr)\frac{\rho
^3_v}{v}\bigg|_{v=0}^{t-s}
- \int_0^{t-s} Y'(v)
\frac{\rho
^3_v}{v}\,\intd v
\\
& =& Y(t)\frac{\rho^3_{t-s}}{t-s} - \int_0^{t-s}
Y'(v)\frac{\rho
^3_v}{v}\,\intd v,
\end{eqnarray*}
where the integral in the last line is understood as a Bochner
integral. Since the total mass of the measure $\rho^3_v$ grows
quadratically in $v$, the last equality tells us that under the
conditions on $Y$, $\caK_g(Y)(t,s)$ is well-defined as a linear
operator on $\caH$. In this case, $\caK_g(Y)(t,s)$ exists as a linear
operator on $\caH$ defined by
\[
\caK_g(Y) (t,s)\sigma(s) = Y(t)\frac{\rho^3_{t-s}(\sigma(s))}{t-s} -
\int
_0^{t-s} Y'(v)\frac{\rho^3_v(\sigma(s))}{v}\,
\intd v,
\]
where
%
%
\begin{equation}
\rho^3_v\bigl(\sigma(s)\bigr) = \int
_{\partial B_3(0,v)} \sigma(s,y)\rho^3_v(\intd y)
= v^2\int_{\partial B_3(0,1)}\sigma(s,vy)\rho
^3_1(\intd y). \label{eqrho}
\end{equation}
Under these conditions the $X$-integral is given by
\begin{eqnarray*}
\int_0^t Y(s)\,\intd X(s) &= & \int
_0^t \biggl(Y(t)\frac{\rho
^3_{t-s}}{t-s} - \int
_0^{t-s} Y'(v)\frac{\rho^3_v}{v}\,
\intd v \biggr)\sigma(s)\delta B(s)
\\
&&{} + \tr_\caH\int_0^t
\biggl(D_sY(t)\frac{\rho^3_{t-s}}{t-s} - \int_0^{t-s}
D_sY'(v)\frac{\rho^3_v}{v}\,\intd v \biggr)\sigma(s)\,\intd s,
\end{eqnarray*}
and this holds for any $\sigma$ for which (\ref{eqrho}) is finite
and for which the integrability conditions from Definition \ref
{defkernel} hold.
\end{pf}

The same procedure can be done in even higher dimensions $d$. The
fundamental solution to the wave equation becomes in these cases
\[
g(t) = \cases{ \displaystyle\frac{2\uppi^{d/2}}{\Gamma(d/2)}1_{ t>0 }
\biggl(
\frac
{1}{t}\frac{\partial}{\partial t} \biggr)^{(d-3)/2}\frac{\rho
_t^d}{t}, &
\quad if $d$ is odd,
\vspace*{3pt}\cr
\displaystyle\frac{2\uppi^{d/2}}{\Gamma(d/2)}1_{ t>0 } \biggl(
\frac
{1}{t}\frac{\partial}{\partial t} \biggr)^{(d-2)/2} \bigl(t^2-
\llvert x\rrvert^2\bigr)_+^{-1/2}, &\quad if $d$ is even.}
\]
This can also be treated by integration by parts. Let us show this for
the case $d=4$, where all the higher-order cases follow in the same
manner. For $d=4$, we deal with the regular distribution $(t^2-\llvert
\cdot
\rrvert ^2)_+^{-1/2}$ given for each test function $f$ by
\[
\bigl(t^2-\llvert\cdot\rrvert^2\bigr)_+^{-1/2}f
= \int_{B_4(0,t)} f(y)\frac{\intd y}{\sqrt
{t^2-\llvert y\rrvert ^2}} = t^2\int
_{B_4(0,1)} f(ty)\frac{\intd y}{\sqrt
{1-\llvert y\rrvert ^2}}.
\]

%
\begin{proposition}
Suppose that Assumption \ref{assdefintegral} and the integrability
conditions for $\sigma$ in Definition \ref{defkernel} hold. Assume
that $Y$ be a random element with values in the linear functionals on
$\caH$ such that $Y$ is differentiable and $Y\in \mathrm{o}(v^{-1})$, $Y'\in
\mathrm{o}(v^{-2})$ and $Y''\in \mathrm{o}(v^{-3})$ at zero. Then, for the special choice
of $g(t,s)=c\rho^4_{t-s}/(t-s)$, $Y\in\caI^X(0,t)$.
\end{proposition}

\begin{pf}
We compute as in the previous proposition
%
%
\begin{eqnarray} \label{eqproofd=4}
&& c_4^{-1}\caK_g(Y) (t,s)
\nonumber
\\
&&\quad = Y(s)g(t-s) + \int_0^{t-s} \bigl(Y(v+s)-Y(s)
\bigr)\frac{\intd
}{\intd v} \biggl(\frac{1}{v}\frac{\intd}{\intd v} \biggr)
\bigl(v^2-\llvert\cdot\rrvert^2\bigr)_+^{-1/2}\,
\intd v
\nonumber\\[-8pt]\\[-8pt]\nonumber
&&\quad = Y(t)g(t-s) - \frac{Y'(t-s)}{t-s} \bigl((t-s)^2-\llvert\cdot
\rrvert^2 \bigr)^{-1/2}_+
\nonumber
\\
&&\qquad {} + \int_0^{t-s}
\frac
{Y''(v)v-Y'(v)}{v^2}\bigl(v-\llvert\cdot\rrvert^2\bigr)_+^{-1/2}\,
\intd v.\nonumber
\end{eqnarray}
In order for this to be well-defined, we need the following condition
for an integration by parts procedure:
\[
\lim_{v\downarrow0} \frac{Y'(v)}{v}\bigl(v^2-\llvert
\cdot\rrvert^2\bigr)_+^{-1/2} = 0.
\]
Under the conditions on $Y$, all the terms in (\ref{eqproofd=4}) are
well-defined.
\end{pf}

Now we want to go a step further and show that even the random-field
integral exists. We show this in the case $d=3$, but the method can be
applied in higher dimensions, also. So we assume in the following that
$Y(s,\ast)$ vanishes at infinity for all $s\in[0,T]$ and that it is
differentiable in $s$ and in $x$. Then we can apply an integration by
parts procedure and obtain
\[
\int_{\R^3} Y(s,z)g(t,s;\intd z,y) = \int
_{\R^3} g(t,s;z,y)\partial^3_zY(s,z)\,\intd z,
\]
where $\partial^3_z=\partial^3/\partial z_1\partial z_2\partial z_3$.
Under these conditions, we can rewrite the integral with respect to the
martingale measure and use its isometry to obtain
%
%
\begin{eqnarray}\label{eqproofrfw}
&& \E\biggl[ \biggl(\int_0^t\int
_{\mathbb{R}^d}\caK_g(Y) (t,s,y)\sigma(s,y)M(\intd s,
\intd y) \biggr)^2 \biggr]
\nonumber\\[-8pt]\\[-8pt]\nonumber
&&\quad = \E\biggl[\int_0^t\int
_{\mathbb{R}^d}\int_{\mathbb{R}^d}\caK_g(Y)
(t,s,y-z)\caK_g(Y) (t,s,y)\sigma(s,y-z)\sigma(s,y)\,\intd y\,\Gamma(
\intd z)\,\intd s \biggr].\quad
\end{eqnarray}
Note that by integration by parts with respect to the temporal and
spatial argument
%
%
\begin{eqnarray} \label{eqproofrfw2}
&& \caK_g(Y) (t,s,y-z)\caK_g(Y) (t,s,y)
\nonumber
\\
&&\quad = \biggl(\int_{\R^3} \partial^3_{x^1}Y
\bigl(t,x^1\bigr)g\bigl(t-s;x^1-y+z\bigr)\,\intd
x^1
\nonumber\\[-8pt]\\[-8pt]\nonumber
&&\hspace*{5pt}\qquad {} + \int_0^{t-s}\int
_{\R^3} \partial^4_{v,x^1}Y
\bigl(v,x^1\bigr)g\bigl(v;\intd x^1-y+z\bigr)\,\intd
x^1\,\intd v \biggr)
\nonumber
\\
&&\qquad {}\times\biggl(\int_{\R^3} \partial
^3_{x^2}Y\bigl(t,x^2\bigr)g
\bigl(t-s;x^2-y\bigr)\,\intd x^2
+ \int_0^{t-s}\int
_{\R^3} \partial^4_{v,x^2}Y
\bigl(v,x^2\bigr)g\bigl(v;x^2-y\bigr)\,\intd x^2\,
\intd v \biggr).\nonumber
\end{eqnarray}
Factoring out this product and substituting this into (\ref
{eqproofrfw}) yields to four stochastic integrals of which we will
explicitly only treat the first one. So we have
\begin{eqnarray*}
&& \E\biggl[\int_0^t\int_{\R^3}
\int_{\R^3} \biggl(\int_{\R
^3}
\partial^3_{x^1}Y\bigl(t,x^1\bigr)g
\bigl(t-s,x^1-y+z\bigr)\,\intd x^1 \biggr)
\nonumber
\\
&&\quad {} \times\biggl(\int_{\R^3}\partial
^3_{x^2}Y\bigl(t,x^2\bigr)g
\bigl(t-s,x^2-y\bigr)\,\intd x^2 \biggr)\sigma(s,y-z)
\sigma(s,y)\,\intd y\,\Gamma(\intd z)\,\intd s \biggr]
\\
&&\quad \leq\int_{\R^3}\int_{\R^3} \sup
_{(r,y)\in[0,T]\times\R^3}\E\bigl[\bigl\llvert\partial^3_{x^1}Y
\bigl(t,x^1\bigr)\partial^3_{x^2}Y
\bigl(t,x^2\bigr)\bigr\rrvert\bigl\llvert\sigma(r,y)\bigr\rrvert
^2 \bigr]
\nonumber
\\
&&\qquad {}\times\int_0^t\int
_{\R^3} \bigl\llvert\exp\bigl(\ii\xi\bigl(x^1+x^2
\bigr)\bigr)\bigr\rrvert\bigl\llvert\caF\Lambda(t-s) (\xi)\bigr\rrvert
^2\mu(\intd\xi)\,\intd s\,\intd x^1\,\intd x^2
\\
&&\quad \leq\int_{\R^3} \bigl(\E\bigl[\bigl\llvert
\partial_x^3 Y(t,x)\bigr\rrvert^{2p} \bigr]
\bigr)^{1/p}\,\intd x
\nonumber
\\
&&\hspace*{12pt}\qquad {}\times\sup_{(r,y)\in[0,T]\times\R^3} \E\bigl[\bigl\llvert
\sigma(r,y)\bigr\rrvert^{2q} \bigr]^{1/q}\int
_0^t\int_{\R^3} \bigl\llvert
\caF\Lambda(t-s) (\xi)\bigr\rrvert^2\mu(\intd\xi)\,\intd s,
\end{eqnarray*}
where in the last line we used H\"older's inequality with
$p^{-1}+q^{-1}=1$. The other three terms in (\ref{eqproofrfw2}) and
the pathwise integral involving the Malliavin derivative can be treated
in a similar way, so that at the end the random field integral $\int
_0^t\int_{\R^3} Y(s,y)X(\intd s,\intd y)$ exists when $g$ is the
fundamental solution to the three-dimensional wave equation if $\sigma
$ has uniform (in time and space) $2q$-moments, if
\[
\int_0^t\int_{\R^3} \bigl
\llvert\caF\Lambda(t-s) (\xi)\bigr\rrvert^2\mu(\intd\xi)\,\intd s <
\infty,
\]
for $\caF$ being the Fourier transform, and if
\[
\int_{\R^3} \bigl(\E\bigl[\bigl\llvert\partial_x^3
Y(t,x)\bigr\rrvert^{2p} \bigr] \bigr)^{1/(2p)} + \int
_0^{t-s}\int_{\R^3} \bigl(\E
\bigl[\bigl\llvert\partial_{u,x}^4 Y(u,x)\bigr\rrvert
^{2p} \bigr] \bigr)^{1/p}\,\intd x\,\intd u < \infty. 
\]
Here, we see a certain trade-off between the integrability conditions
of $Y$ and $\sigma$, if we let one of the two to be more irregular,
the other one has to be more regular. A similar condition with the
Malliavin derivatives of $Y$ has to hold so that the pathwise integral
is well-efined. All this means that $Y$ has to be in some Sobolev-type
space where the (weak) derivatives of the process and of the Malliavin
derivatives of the process exist and an integrability condition has to
hold. These conditions are far from optimal, but serve as a first step
toward a general existence theorem for the random-field $X$-integral
when $g$ is a distribution.


\section*{Acknowledgements}
The authors would like to thank Marta Sanz-Sol\'e for reading a
preliminary version of this article.

Fred Espen Benth is grateful for the financial support from the projects
\textit{Managing weather risk in electricity markets} (\textit
{MAWREM}) and
\textit{Energy markets}: \textit{Modelling}, \textit{optimization
and simulation} (\textit{EMMOS}),
both sponsored by the Norwegian Research Council under the RENERGI and
EVITA programmes, respectively.

Andr{\'e} S{\"u}\ss{} is supported by the grant MICINN-FEDER MTM
2009-07203 from the Direcci\'on
General de Investigaci\'on, Ministerio de Economia y Competitividad, Spain.
The research leading to this article has been done while the second
author was on a {\it Estanc\'\i{}a breve} program at CMA in Oslo funded
by the Spanish Department of Science under EEBB-I-12-04650.
The financial support and the hospitality by the CMA staff is
gratefully acknowledged.


%

\printhistory

\begin{thebibliography}{19}

\bibitem{alosnualart}
%
\begin{barticle}[mr]
\bauthor{\bsnm{Al{\`o}s},~\bfnm{Elisa}\binits{E.}},
\bauthor{\bsnm{Mazet},~\bfnm{Olivier}\binits{O.}} \AND
\bauthor{\bsnm{Nualart},~\bfnm{David}\binits{D.}}
(\byear{2001}).
\btitle{Stochastic calculus with respect to {G}aussian processes}.
\bjournal{Ann. Probab.}
\bvolume{29}
\bpages{766--801}.
\bid{doi={10.1214/aop/1008956692}, issn={0091-1798}, mr={1849177}}
\end{barticle}
%

\bptok{imsref}%
\endbibitem

\bibitem{vmlv}
%
\begin{barticle}[mr]
\bauthor{\bsnm{Barndorff-Nielsen},~\bfnm{Ole~E.}\binits{O.E.}},
\bauthor{\bsnm{Benth},~\bfnm{Fred~Espen}\binits{F.E.}},
\bauthor{\bsnm{Pedersen},~\bfnm{Jan}\binits{J.}} \AND
\bauthor{\bsnm{Veraart},~\bfnm{Almut~E.~D.}\binits{A.E.D.}}
(\byear{2014}).
\btitle{On stochastic integration for volatility modulated L\'
evy-driven {V}olterra processes}.
\bjournal{Stochastic Process. Appl.}
\bvolume{124}
\bpages{812--847}.
\bid{doi={10.1016/j.spa.2013.09.007}, issn={0304-4149}, mr={3131315}}
\end{barticle}
%

\bptok{imsref}%
\endbibitem

\bibitem{benth}
%
\begin{bincollection}[mr]
\bauthor{\bsnm{Barndorff-Nielsen},~\bfnm{Ole~E.}\binits{O.E.}},
\bauthor{\bsnm{Benth},~\bfnm{Fred~Espen}\binits{F.E.}} \AND
\bauthor{\bsnm{Veraart},~\bfnm{Almut~E.~D.}\binits{A.E.D.}}
(\byear{2011}).
\btitle{Ambit processes and stochastic partial differential equations}.
In \bbooktitle{Advanced Mathematical Methods for Finance}
\bpages{35--74}.
\blocation{Heidelberg}:
\bpublisher{Springer}.
\bid{doi={10.1007/978-3-642-18412-3_2}, mr={2752540}}
\bptnote{check pages}%
\end{bincollection}
%

\bptok{imsref}%
\endbibitem

\bibitem{carmonatehranchi}
%
\begin{bbook}[mr]
\bauthor{\bsnm{Carmona},~\bfnm{Ren{\'e}~A.}\binits{R.A.}} \AND
\bauthor{\bsnm{Tehranchi},~\bfnm{Michael~R.}\binits{M.R.}}
(\byear{2006}).
\btitle{Interest Rate Models: An Infinite Dimensional Stochastic
Analysis Perspective}.
\bseries{Springer Finance}.
\blocation{Berlin}:
\bpublisher{Springer}.
\bid{mr={2235463}}
\end{bbook}
%

\bptok{imsref}%
\endbibitem

\bibitem{dalang}
%
\begin{barticle}[mr]
\bauthor{\bsnm{Dalang},~\bfnm{Robert~C.}\binits{R.C.}}
(\byear{1999}).
\btitle{Extending the martingale measure stochastic integral with
applications to spatially homogeneous s.p.d.e.'s}.
\bjournal{Electron. J. Probab.}
\bvolume{4}
\bpages{1--29 (electronic)}.
\bid{doi={10.1214/EJP.v4-43}, issn={1083-6489}, mr={1684157}}
\bptnote{check pages}%
\end{barticle}
%

\bptok{imsref}%
\endbibitem

\bibitem{dalangquer}
%
\begin{barticle}[mr]
\bauthor{\bsnm{Dalang},~\bfnm{Robert~C.}\binits{R.C.}} \AND
\bauthor{\bsnm{Quer-Sardanyons},~\bfnm{Llu{\'{\i}}s}\binits{L.}}
(\byear{2011}).
\btitle{Stochastic integrals for spde's: A comparison}.
\bjournal{Expo. Math.}
\bvolume{29}
\bpages{67--109}.
\bid{doi={10.1016/j.exmath.2010.09.005}, issn={0723-0869}, mr={2785545}}
\end{barticle}
%

\bptok{imsref}%
\endbibitem

\bibitem{dapratozabczyk}
%
\begin{bbook}[mr]
\bauthor{\bsnm{Da Prato},~\bfnm{Giuseppe}\binits{G.}} \AND
\bauthor{\bsnm{Zabczyk},~\bfnm{Jerzy}\binits{J.}}
(\byear{1992}).
\btitle{Stochastic Equations in Infinite Dimensions}.
\bseries{Encyclopedia of Mathematics and Its Applications}
\bvolume{44}.
\blocation{Cambridge}:
\bpublisher{Cambridge Univ. Press}.
\bid{doi={10.1017/CBO9780511666223}, mr={1207136}}
\bptnote{check year}%
\end{bbook}
%

\bptok{imsref}%
\endbibitem

\bibitem{diesteluhl}
%
\begin{bbook}[mr]
\bauthor{\bsnm{Diestel},~\bfnm{J.}\binits{J.}} \AND
\bauthor{\bsnm{Uhl},~\bfnm{J.~J.}\binits{J.J.} \bsuffix{Jr.}}
(\byear{1977}).
\btitle{Vector Measures}.
\blocation{Providence, RI}:
\bpublisher{Amer. Math. Soc.}
\bid{mr={0453964}}
\end{bbook}
%

\bptok{imsref}%
\endbibitem

\bibitem{dunfordschwartz}
%
\begin{bbook}[mr]
\bauthor{\bsnm{Dunford},~\bfnm{Nelson}\binits{N.}} \AND
\bauthor{\bsnm{Schwartz},~\bfnm{Jacob~T.}\binits{J.T.}}
(\byear{1988}).
\btitle{Linear Operators. {P}art I: General Theory}.
\bseries{Wiley Classics Library}.
\blocation{New York}:
\bpublisher{Wiley}.
\bid{mr={1009162}}
\end{bbook}
%

\bptok{imsref}%
\endbibitem

\bibitem{engelnagel}
%
\begin{bbook}[mr]
\bauthor{\bsnm{Engel},~\bfnm{Klaus-Jochen}\binits{K.-J.}} \AND
\bauthor{\bsnm{Nagel},~\bfnm{Rainer}\binits{R.}}
(\byear{2006}).
\btitle{A Short Course on Operator Semigroups}.
\bseries{Universitext}.
\blocation{New York}:
\bpublisher{Springer}.
\bid{mr={2229872}}
\end{bbook}
%

\bptok{imsref}%
\endbibitem

\bibitem{grorudpardoux}
%
\begin{barticle}[mr]
\bauthor{\bsnm{Grorud},~\bfnm{A.}\binits{A.}} \AND
\bauthor{\bsnm{Pardoux},~\bfnm{{\'E}.}\binits{{\'E}.}}
(\byear{1992}).
\btitle{Int\'egrales hilbertiennes anticipantes par rapport \`a un
processus de {W}iener cylindrique et calcul stochastique associ\'e}.
\bjournal{Appl. Math. Optim.}
\bvolume{25}
\bpages{31--49}.
\bid{doi={10.1007/BF01184155}, issn={0095-4616}, mr={1133251}}
\end{barticle}
%

\bptok{imsref}%
\endbibitem

\bibitem{sergelang}
%
\begin{bbook}[mr]
\bauthor{\bsnm{Lang},~\bfnm{Serge}\binits{S.}}
(\byear{1993}).
\btitle{Real and Functional Analysis},
\bedition{3rd} ed.
\bseries{Graduate Texts in Mathematics}
\bvolume{142}.
\blocation{New York}:
\bpublisher{Springer}.
\bid{doi={10.1007/978-1-4612-0897-6}, mr={1216137}}
\end{bbook}
%

\bptok{imsref}%
\endbibitem

\bibitem{metivier}
%
\begin{bbook}[mr]
\bauthor{\bsnm{M{\'e}tivier},~\bfnm{Michel}\binits{M.}}
(\byear{1982}).
\btitle{Semimartingales: A Course on Stochastic Processes}.
\bseries{De Gruyter Studies in Mathematics}
\bvolume{2}.
\blocation{Berlin}:
\bpublisher{de Gruyter}.
\bid{mr={0688144}}
\end{bbook}
%

\bptok{imsref}%
\endbibitem

\bibitem{nualart}
%
\begin{bbook}[mr]
\bauthor{\bsnm{Nualart},~\bfnm{David}\binits{D.}}
(\byear{2006}).
\btitle{The {M}alliavin Calculus and Related Topics},
\bedition{2nd} ed.
\bseries{Probability and Its Applications (New York)}.
\blocation{Berlin}:
\bpublisher{Springer}.
\bid{mr={2200233}}
\end{bbook}
%

\bptok{imsref}%
\endbibitem

\bibitem{nualartpardoux}
%
\begin{barticle}[mr]
\bauthor{\bsnm{Nualart},~\bfnm{D.}\binits{D.}} \AND
\bauthor{\bsnm{Pardoux},~\bfnm{{\'E}.}\binits{{\'E}.}}
(\byear{1988}).
\btitle{Stochastic calculus with anticipating integrands}.
\bjournal{Probab. Theory Related Fields}
\bvolume{78}
\bpages{535--581}.
\bid{doi={10.1007/BF00353876}, issn={0178-8051}, mr={0950346}}
\end{barticle}
%

\bptok{imsref}%
\endbibitem

\bibitem{pronkveraar}
%
\begin{barticle}[mr]
\bauthor{\bsnm{Pronk},~\bfnm{Matthijs}\binits{M.}} \AND
\bauthor{\bsnm{Veraar},~\bfnm{Mark}\binits{M.}}
(\byear{2014}).
\btitle{Tools for {M}alliavin calculus in UMD {B}anach spaces}.
\bjournal{Potential Anal.}
\bvolume{40}
\bpages{307--344}.
\bid{doi={10.1007/s11118-013-9350-0}, issn={0926-2601}, mr={3201985}}
\bptnote{check volume, check pages, check year}%
\end{barticle}
%

\bptok{imsref}%
\endbibitem

\bibitem{protter}
%
\begin{bbook}[mr]
\bauthor{\bsnm{Protter},~\bfnm{Philip~E.}\binits{P.E.}}
(\byear{2005}).
\btitle{Stochastic Integration and Differential Equations}.
\bseries{Stochastic Modelling and Applied Probability}
\bvolume{21}.
\blocation{Berlin}:
\bpublisher{Springer}.
\bid{doi={10.1007/978-3-662-10061-5}, mr={2273672}}
\end{bbook}
%

\bptok{imsref}%
\endbibitem

\bibitem{sanzsuess1}
%
\begin{barticle}[mr]
\bauthor{\bsnm{Sanz-Sol{\'e}},~\bfnm{Marta}\binits{M.}} \AND
\bauthor{\bsnm{S{\"u}{\ss}},~\bfnm{Andr{\'e}}\binits{A.}}
(\byear{2013}).
\btitle{The stochastic wave equation in high dimensions: {M}alliavin
differentiability and absolute continuity}.
\bjournal{Electron. J. Probab.}
\bvolume{18}
\bpages{no. 64, 28}.
\bid{doi={10.1214/EJP.v18-2341}, issn={1083-6489}, mr={3078023}}
\bptnote{check pages}%
\end{barticle}
%

\bptok{imsref}%
\endbibitem

\bibitem{walsh}
%
\begin{bincollection}[mr]
\bauthor{\bsnm{Walsh},~\bfnm{John~B.}\binits{J.B.}}
(\byear{1986}).
\btitle{An introduction to stochastic partial differential equations}.
In \bbooktitle{\'Ecole D'\'et\'e de Probabilit\'es de {S}aint-{F}lour,
XIV~-- 1984}.
\bseries{Lecture Notes in Math.}
\bvolume{1180}
\bpages{265--439}.
\blocation{Berlin}:
\bpublisher{Springer}.
\bid{doi={10.1007/BFb0074920}, mr={0876085}}
\bptnote{check pages}%
\end{bincollection}
%

\bptok{imsref}%
\endbibitem

\end{thebibliography}
\end{document}